\documentclass [reqno,12pt,a4paper]{amsart} 

\usepackage{graphicx,color,amsmath,amssymb}
\usepackage{subfigure}
\usepackage{natbib}
\begin{document}

\def \RR {{\mathbb R}}
\def \ZZ {{\mathbb Z}}
\def \NN {{\mathbb N}}
\def \vx {{^x\!v}}
\def \vy {{^y\!v}}
\def \dx {{\Delta X}}
\def \dy {{\Delta Y}}
\def \dt {{\Delta t}}
\def \dtc {{\Delta t_c}}
\def \dtcfl {{\Delta t_c}}
\def \ra {{\rightarrow}}
\def \xip {{x^\kappa_{j+1/2,l}}}
\def \xfp {{x^\kappa_{j+1/2,r}}}
\def \xim {{x^\kappa_{j-1/2,l}}}
\def \xfm {{x^\kappa_{j-1/2,r}}}
\def \yip {{y^\kappa_{k+1/2,l}}}
\def \yfp {{y^\kappa_{k+1/2,r}}}
\def \yim {{y^\kappa_{k-1/2,l}}}
\def \yfm {{y^\kappa_{k-1/2,r}}}
\def \tmk {{t^\kappa_m}}
\def \cal {{\mathcal}}
\def \II {{\mbox{I\!I}}}
\def \III {{\mbox{I\!I\!I}}}
\def \IX {{{I\hspace{-0.09cm}X}}}
\def \SD {{SD1$-$2D}}
\def \SDD {{SD2$-$2D}}

\title{A New two-dimensional Second Order
  Non-oscillatory Central Scheme Applied to multiphase flows
  in heterogeneous porous media} \author{F. Furtado$^{2}$,
  F. Pereira$^{1}$ and S. Ribeiro$^{1}$} \date{}
\maketitle 

\noindent {$^{1}$ State University of Rio de Janeiro, Nova
  Friburgo, RJ 28630-050, Brazil.}

\noindent {$^{2}$ University of Wyoming, Laramie, WY
  82071-3036, U.S.A.}

\pagestyle{myheadings}\markboth{F. Furtado, F. Pereira and
S. Ribeiro}{SD2-2D}

\section*{Abstract}

We are concerned with central differencing schemes for
solving scalar hyperbolic conservation laws.
We compare the Kurganov-Tadmor (KT) two-dimensional
\citep{KT2000} second order semi-discrete central scheme
in {\em dimension by dimension} formulation with a new {\em
  two-dimensional} approach introduced here and
applied in numerical simulations for two-phase,
two-dimensional flows in heterogeneous formations.  
This semi-discrete central scheme is based on the
ideas of Rusanov's method (\cite{RUSANOV}) using a more
precise information about the local speeds of wave propagation
computed at each Riemann Problem in two-space dimensions.
We find
the KT dimension by dimension has a much simpler
mathematical description than the genuinely two-dimensional
one with a little more numerical diffusion, particularly in the
presence of viscous fingers. Unfortunately, as one can see
in \cite{AFPR2006}, the KT with the dimension by dimension
approach might produce incorrect boundary behavior in a
typical geometry used in the study of porous media flows:
the quarter of a five spot. This problem has been corrected
by the authors with the new semi-discrete scheme proposed
here.   We conclude with numerical examples of
two-dimensional, two-phase flow associated with two distinct
flooding problems: a
two-dimensional flow in a rectangular heterogeneous
reservoir (called slab geometry)  and a
two-dimensional flow in a 5-spot geometry homogeneous
reservoir.


\section{INTRODUCTION}\label{sec1}

We consider a model for two-phase flow, immiscible and
incompressible displacement in heterogeneous porous
media. The highly nonlinear governing equations are of very
practical importance in petroleum engineering
\citep{dwpb,CJ1986} (see also \citep{fpcross} and the
references therein for recent studies for the scale-up
problem for such equations).

The conventional theoretical description of two-phase flow
in a porous medium, in the limit of vanishing capillary
pressure, is via Darcy's law coupled to the Buckley-Leverett
equation. The two phases will be referred to as water and
oil, and indicated by the subscripts $w$ and $o$,
respectively.  We also assume that the two fluid phases
saturate the pores. With no sources or sinks, and neglecting
the effects of gravity these equations are (see Peaceman
\citep{dwpb}):
\begin{equation}
\nabla \cdot {\bf v}=0, \quad {\bf v}=-\lambda(s) K({\bf x}) \nabla p,  
\label{eq:preeq}
\end{equation} 
\begin{equation}
\frac{\partial s}{\partial t}+ \nabla \cdot (f(s){\bf v})=0.  
\label{eq:sateq}
\end{equation} 
Here, ${\bf v}$ is the total seepage velocity, $s$ is the
water saturation, $K({\bf x})$ is the absolute permeability,
and $p$ is the pressure. The constant porosity has been
scaled out by a change of the time variable. The
constitutive functions $\lambda(s)$ and $f(s)$ represent the
total mobility and the fractional flow of water,
respectively.

We are concerned with numerical schemes for solving scalar
hyperbolic conservation laws arising in the simulation of
multiphase flows in multidimensional heterogeneous porous
media. These schemes are non-oscillatory and enjoy the main
advantage of Godunov-type central schemes: simplicity, i.e.,
they employ neither characteristic decomposition nor
approximate Riemann solvers. This makes them universal
methods that can be applied to a wide variety of physical
problems, including hyperbolic systems of conservation laws.
The two main classes of Godunov methods are upwind and
central schemes.

The Lax-Friedrichs (LxF) scheme \citep{LAX1954} is the
canonical first order central scheme, which is the
forerunner of all central differencing schemes. It is a
non-oscillatory scheme based on piecewise constant
approximate solution and it also enjoys simplicity.
Unfortunately the excessive numerical dissipation in the LxF
recipe (of order ${\mathcal O}(\Delta X^{2}/\Delta t)$)
yields a poor resolution, which seems to have delayed the
development of high resolution central schemes when compared
with the earlier developments of the high resolution upwind
methods.  Only in 1990 a second order generalization to the
LxF scheme was introduced in \citep{TN}. They used a
staggered form of the LxF scheme and replaced the first
order piecewise constant solution with a van Leer's
MUSCL-type piecewise linear second order approximation
\citep{VANLEER1979}.  The numerical dissipation present in
this new central scheme has an amplitude of order ${\mathcal
  O}({\Delta X}^{4}/\Delta t)$. (see \citep{AFF},
\citep{AFMPb}, \citep{AFMP}, \citep{ADFMP},
for recent studies in three phase flows using the Nessyahu
and Tadmor (NT) central scheme). In spite of the good
resolution obtained by the Nessyahu and Tadmor scheme, much
higher than in the first order LxF scheme, there are still
some difficulties with small time steps which arise, e.g. in
multiphase flows in highly heterogeneous petroleum
reservoirs or aquifers. To overcome this difficulty, we can
use a semi-discrete formulation coupled with an appropriate
ODE solver retaining simplicity and high resolution with
lower numerical viscosity, proportional to the vanishing
size of the time step $\Delta t$. Both LxF and NT schemes do
not admit a semi-discrete form; see \citep{KT2000} for a
detailed description of the one-dimensional Kurganov and
Tadmor central scheme which is the first fully discrete
Godunov-Type central scheme admitting a semi-discrete form.

We compare the Kurganov-Tadmor (KT) two-dimensional
\citep{KT2000} second order semi-discrete central scheme
in {\em dimension by dimension} formulation with a {\em genuinely
two-dimensional} approach applied in numerical simulations
for two-phase, two-dimensional flows in heterogeneous
formations.  We find the KT dimension by dimension has a
much simpler mathematical description than the genuinely
two-dimensional one adding only a little more diffusion,
particularly in the presence of viscous
fingers. Unfortunately, the KT with the dimension by
dimension approach might produce incorrect boundary behavior
in a typical geometry used in the study of porous media
flows: the quarter of a five spot. These results are
presented in \citep{AFPR2006}.  This problem motivated the
authors to develop a genuinely two-dimensional formulation
which is then presented in section
(\ref{central_scheme}). Although a similar two-dimensional
formulation was available in a early work \citep{MR1836876},
ours was developed independently to deal with two-phase
flows, immiscible and incompressible displacement in
heterogeneous porous media.  It shares the same
general ideas with the work of Kurganov-Petrovna but differs
in many technical details.

This paper is organized as follows. In Section \ref{model}
we introduce the model for two-phase flows, immiscible and
incompressible displacement in heterogeneous porous
media. In Section \ref{central_scheme} we discuss the
mathematical formulation for the KT central scheme in
dimension by dimension approach and in a genuinely
two-dimensional one.  In Section \ref{sec:aplication} we will
present some numerical results when we apply the KT central
differencing scheme with both approaches mentioned above to
porous media flows.

\section{NUMERICAL APPROXIMATION OF TWO-PHASE FLOW}\label{model}

\subsection{Operator splitting for two-phase flow.}
An operator splitting technique is employed for the
computational solution of the saturation equation
(\ref{eq:sateq}) and the pressure equation \eqref{eq:preeq} which
are solved sequentially with distinct time steps. This
method has proved to be computationally efficient in
producing accurate numerical solutions for two-phase flow
\citep{dfp}.

Typically, for computational efficiency, larger time steps
are used to compute the pressure \eqref{eq:preeq}. The
splitting allows time steps used in the solution of the
pressure-velocity equation that are longer than those
allowed in the convection calculation {~\eqref{eq:sateq}}.
Thus, we introduce two time steps: $\Delta t_c$ used to
compute the solution of the hyperbolic problem, and $\Delta
t_p$ used in the pressure-velocity calculation such that $\Delta
t_p \geq \Delta t_c$. We remark that in practice, variable
time steps are always useful, especially for the convection
micro-steps subject dynamically to a $CFL$ condition.

For the global pressure solution (the pertinent elliptic
equation), we use a (locally conservative) hybridized mixed
finite element discretization equivalent to cell-centered
finite differences \citep{dfp}, which
effectively treats the rapidly changing permeabilities that
arise from stochastic geology and produces accurate velocity
fields. The pressure and the Darcy velocity are approximated
at times $t^m = m\Delta t_p$, $m=0,1,2,\dots$. The algebraic
system resulting from the discretized equations can be
solved by a preconditioned conjugate gradient procedure
(PCG) or by a multi-grid procedure (\citep{dfp}).

We use high resolution numerical central scheme (see
\citep{KT2000}) for solving the scalar hyperbolic
conservation laws arising in the convection of the fluid
phases in heterogeneous porous media for two-phase flows -
we will discuss the application of these schemes for
two-phase flows in Section \ref{sec:aplication}.
Theses methods can accurately resolve sharp fronts in the
fluid saturations without introducing spurious oscillations
or excessive numerical diffusion.

The saturation equation is approximated at times
$t^{\kappa}_{m}=t^{m} + {\kappa}\Delta t_c$ for $t^{m} <
\tmk \leq t^{m+1}$ that take into account the
advective transport of water. We will write $t^\kappa$ to
indicate the time step $\tmk$ and $t^{\kappa+1}$ to indicate
$\tmk + \Delta t_c$.


A detailed description of the numerical method that we
employ for the solution of Eqs. (\ref{eq:preeq})-(\ref{eq:sateq})
can be found in \cite{dfp} and references therein (see
also \cite{ADFMP,AFF} and \cite{abreu_splitting_2008} for applications of the
operator splitting technique for three phase flows that
takes into account capillary pressure (diffusive effects)).

{\emph Remark: } To simplify notation, we denote:
\begin{itemize}
\item NT1d for one-dimensional NT scheme;
\item NT2d for two-dimensional NT scheme;
\item KT1d for one-dimensional KT scheme;
\item KTdxd for the KT scheme with dimension by dimension approach and
\item SD2-2D for our two-dimensional approach. 
\end{itemize}

\subsection{TWO SPATIAL DIMENSIONS SECOND ORDER SEMI-DISCRETE CENTRAL SCHEME}
\label{central_scheme}

In this section, we will develop a two-spatial dimension second order 
semi-discrete central scheme (SD2-2D) based on the ideas of \citet{LAX1954,RUSANOV,TN, KT2000} and \citet{NT2D}  which are then applied in numerical approximation of the scalar hyperbolic conservation law
modeling the convective transport of the fluid phases in
two-phase flows and its coupling with lowest order
Raviart-Thomas \citep{prta} locally conservative mixed
finite elements for the associated elliptic
(velocity-pressure part) problem (See \cite{prta}).
We summarize below the basic ideas of the construction of SD2$-$2D numerical scheme:
\begin{itemize}
\item The Lax-Friedrichs method in two-spatial dimensions LxF2D written in the \texttt{REA} algorithm setup  (See \citet{NT2D}) will be used to obtain the two dimensional Rusanov's method SD1-2D. We follow the same procedures presented in \citet{KT2000} in one spatial dimension.
\item The new SD2-2D  numerical scheme will then be obtained as a second order generalization of the SD1-2D. This second order precision is achieved approximating the solution with piecewise linear functions. 
\end{itemize}

\subsection{The staggered non-uniform grid of the SD2-2D central scheme}

\begin{figure}[h]
\centerline{\includegraphics[scale=0.6]{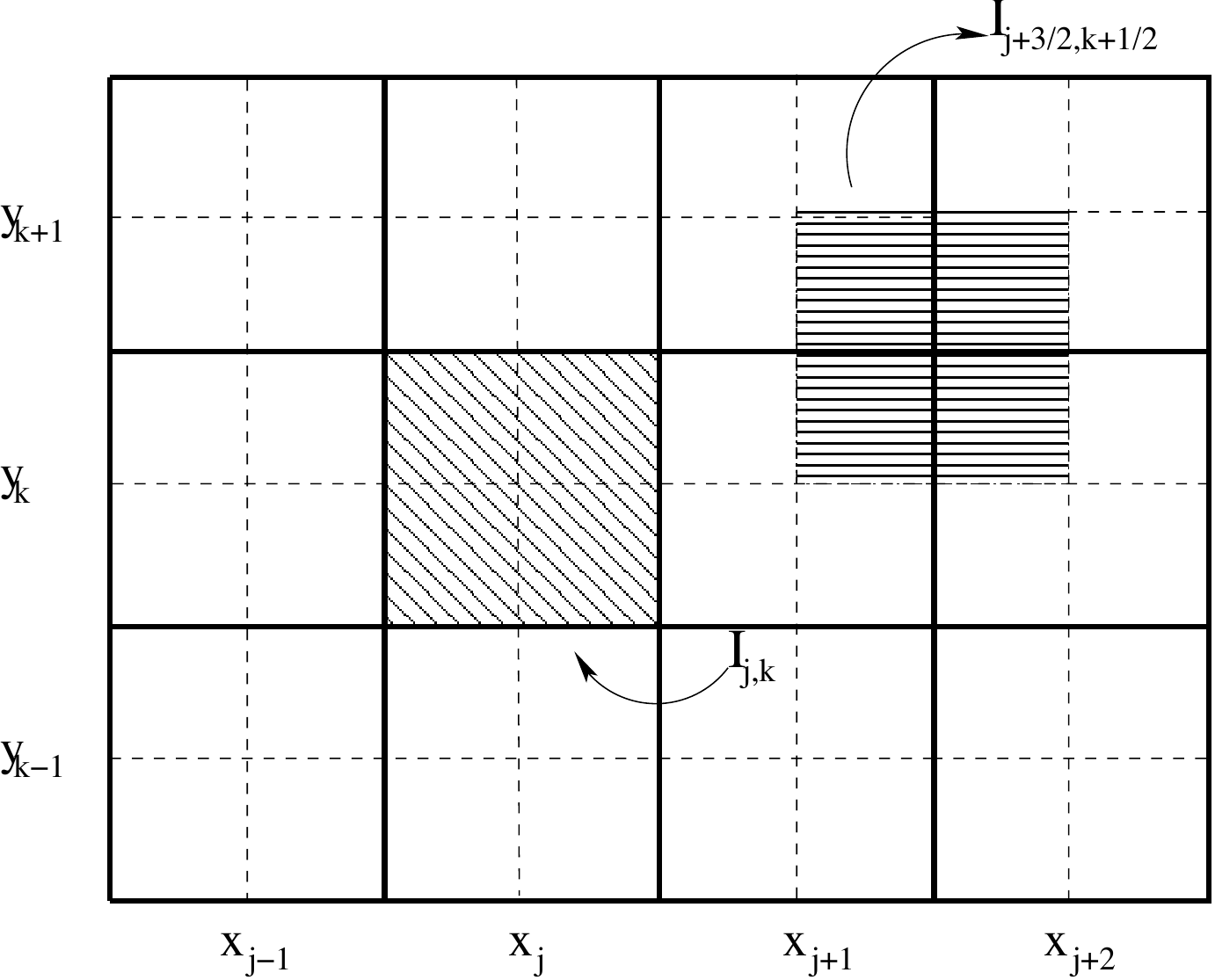}}
\caption{\label{fig:malha_deslocada} The LxF2D uniform grid. The cells  $I_{j,k}$ of the original non-staggered grid are limited by the solid lines and the cells $I_{j+3/2,k+1/2}$ of the staggered grid are limited by the dashed lines.}
\end{figure}

We begin then extending the LxF2D to obtain the SD1-2D following
the same procedures for one dimensional problems.  First, we define the non-staggered and the staggered grids of retangular cells used in the LxF2D. 
The points $(x_j, y_k,t^\kappa)$ are defined as follows.
\begin{eqnarray}
  x_j & = & j\cdot \Delta x, \qquad j=\ldots, -1, 0, 1, \ldots
  \nonumber \\
  y_k & = & k\cdot \Delta y, \qquad k=\ldots, -1, 0, 1, \ldots
  \nonumber \\
  t^\kappa & = & t^m + \kappa\cdot \Delta t_{c}, \qquad
  \kappa=0,\ldots, i, \nonumber
\end{eqnarray}
We denote the cells of the non-staggered grid by $I_{j,k} :=
(x_{j-1/2}, x_{j+1/2})\times (y_{k-1/2}, y_{k+1/2})$. Its area
 $A_{j,k} \equiv \Delta x \cdot \Delta y = (x_{j+1/2} -
x_{j-1/2})\cdot (y_{k+1/2} -y_{k-1/2})$. The time step of the
convective equation \eqref{eq:sateq} is
$\Delta t_{c} = t^{\kappa+1} - t^\kappa$.  The staggered grid is obtained
moving the cells $\Delta x/2$ to the right and $\Delta y/2$ upward
These staggered cells will be denoted by
$I_{j+1/2,k+1/2} := (x_j,
x_{j+1})\times (y_k, y_{k+1})$. Its center is the point
$(x_{j+1/2},y_{k+1/2})$, where  $x_{j+1/2} = x_j + \Delta x/2$ and
$y_{k+1/2} = y_j + \Delta y/2$. The Figure
\ref{fig:malha_deslocada} illustrates  the non-staggered and the staggered grids
showing the cells $I_{j,k}$ e $I_{j+3/2,k+1/2}$ as examples of  non-staggered and  staggered cells, respectively.
 
The scalar hyperbolic conservation law \eqref{eq:sateq} can be written
as
\begin{equation}
\frac{\partial s}{\partial t} + \frac{\partial }{\partial
  x}(\vx f(s)) +
\frac{\partial }{\partial y}(\vy f(s)) = 0,
\label{two2D}
\end{equation}
where $\vx \equiv \vx(x,y,t)$ and $\vy \equiv \vy(x,y,t)$ denote the $x$ and 
$y$ components of the velocity field ${\bf v}$. The
cell averages at time $t^{\kappa}$ are 
\begin{equation}
\displaystyle \overline{S}^{\kappa}_{j,k} := \overline{S}_{j,k}(t^{\kappa})  \equiv \frac{1}{\Delta X\! \Delta Y}
\displaystyle\int_{x_{j-\frac{1}{2}}}^{x_{j+\frac{1}{2}}} 
\displaystyle\int_{y_{k-\frac{1}{2}}}^{y_{k+\frac{1}{2}}} s(x,y,t^{\kappa})\,dxdy.
\label{eq:cellaver}
\end{equation}

The solution $s(x,y,t^\kappa)$ of the problem
\eqref{eq:sateq} at time
$t^\kappa$ is approximated using  piecewise-linear MUSCL-type interpolants (See \citet{VANLEER1979}). 
\begin{equation}
  \label{eq:bilinear_reconstruction}
  \widetilde{S}^\kappa_{j,k}(x,y) =   \overline{S}^\kappa_{j,k} +
  (S_x)^\kappa_{j,k}\cdot (x-x_j) + (S_y)^\kappa_{j,k}\cdot (y-y_k), 
\end{equation}
where $x_{j-1/2} \le x \le x_{j+1/2}$ e $y_{k-1/2} \le y \le
y_{k+1/2}$.  The second-order resolution is  guaranteed if the so-called vector of numerical derivative $(S_x)^\kappa_{j,k}$ and
$(S_y)^\kappa_{j,k}$ satisfy 
\begin{subequations}
\label{eq:dev_2ordem_2D}
\begin{eqnarray}
  \label{eq:dev_2ordem_2D_x_presente}
  (S_x)^\kappa_{j,k} & = & \left.\frac{\partial s}{\partial
    x}\,\right|_{x = x_j,y = y_k,t =t^\kappa} + \mathcal{O}(\Delta X); \\
  \label{eq:dev_2ordem_2D_y_presente}
(S_y)^\kappa_{j,k} & = & \left.\frac{\partial s}{\partial
    y}\right|_{x= x_j,y=y_k,t=t^\kappa} + \mathcal{O}(\Delta Y),   
\end{eqnarray}
\end{subequations}
These numerical derivatives are computed using the MinMod limiter 
\begin{subequations}
 \label{eq:minmod_2D}
\begin{eqnarray}
  \label{eq:minmod_2Dx}
  (S_x)^\kappa_{j,k} & = & \text{MM}\theta\frac{1}{\Delta x}\left\{\overline{S}^\kappa_{j-1,k},\overline{S}^\kappa_{j,k},
    \overline{S}^\kappa_{j+1,k}\right\} \nonumber \\
  & := & \text{MM}\left(\theta \frac{\Delta S^\kappa_{j+1/2,k}}{\Delta x}, 
  \frac{\Delta S^\kappa_{j-1/2,k} -\Delta S^\kappa_{j+1/2,k}}{2\Delta
    x}, \theta \frac{\Delta S^\kappa_{j-1/2,k}}{\Delta x}\right);
\\
 \label{eq:minmod_2Dy}
(S_y)^\kappa_{j,k} & = & \text{MM}\theta\frac{1}{\Delta y}\left\{\overline{S}^\kappa_{j,k-1},\overline{S}^\kappa_{j,k},
    \overline{S}^\kappa_{j,k+1}\right\} \nonumber \\
  & := & \text{MM}\left(\theta \frac{\Delta S^\kappa_{j,k+1/2}}{\Delta y}, 
  \frac{\Delta S^\kappa_{j,k-1/2} -\Delta S^\kappa_{j,k+1/2}}{2\Delta
    y}, \theta \frac{\Delta S^\kappa_{j,k-1/2}}{\Delta y}\right).
\end{eqnarray}
\end{subequations}
Here $\Delta$ denotes the centered difference, $\Delta
S^\kappa_{j+1/2,k} = \overline{S}^\kappa_{j+1,k} - \overline{S}^\kappa_{j,k}$
and the paramter $\theta \in [1,1.8]$
has been chosen in the optimal way in every numerical example with $\theta = 1.8$ beeing the less dissipative limiter. 
The minmod limiter \eqref{eq:minmod_2D} guarantees
the nonoscillatory property and the second-order accuracy.
The reconstruction given by Equations \eqref{eq:bilinear_reconstruction}-\eqref{eq:minmod_2D} also retains conversation, i.e.,
\begin{equation}
  \label{eq:tilde_preserva_conserv}
  \int_{I_{j,k}} \widetilde{S}^\kappa_{j,k}(x,y)\, dx \,dy = \overline{S}^\kappa_{j,k}. 
\end{equation}

\textit{Remark: } We notice that if $(S_x)^\kappa_{j,k}$ and $(S_y)^\kappa_{j,k}$ are equal to zero, then we will get the first-order two-dimensional semi-discrete scheme SD1-2D. Otherwise, we will obtain the second-order two spatial dimensions semi-discrete central scheme SD2-2D.

We consider the model of hyperbolic conservation laws given by
Equation \eqref{eq:sateq} with cell averages as in \eqref{eq:cellaver} and the two-dimensional 
piecewise linear reconstruction defined in 
\eqref{eq:bilinear_reconstruction} and \eqref{eq:dev_2ordem_2D}
with the conservative property \eqref{eq:tilde_preserva_conserv}.
Our goal is to compute an approximated solution
$S_{j,k}(\tmk + \Delta t_c)$ in the original grid
at the next time step. To this end, we apply the Godunov's
algorithm REA. To solve this problem, we integrate the 
conservation law over some control volumes that we need to specify.

{\bf Constructing the staggered nonuniform grid:} 
Kurganov and Tadmor developed 
the KT1D scheme along the lines 
of NT1D (See \citep{KT2000}). The
nonuniform staggered grid in the KT1D was constructed
directly from the 
staggered uniform grid of NT1D with additional information
on the local speeds of wave propagation. In a similar way the 
{\bf nonuniform} staggered grid in the SD2-2D scheme is defined from
the {\bf uniform} staggered grid of the NT2D scheme as follows.

$(1)$ We set $C_{j,k} = [x_{j-1/2}, x_{j+1/2}]
\times [y_{k-1/2},y_{k+1/2}]$ to denote the cells 
of the original non-staggered grid; $C_{j+1/2,k+1/2} = 
[x_j, x_{j+1}] \times [y_k,y_{k+1}]$ to denote the cells
of the uniform staggered grid. We start with 
a piecewise constant approximated solution 
$\bar{S}^\kappa_{j,k}$   
over the original cells $C_{j,k}$.

$(2)$ Next, we move to the staggered uniform grid with cells
given by $C_{j+1/2, k+1/2}$.

\begin{figure}[h]
\centerline{\includegraphics[scale=0.7]{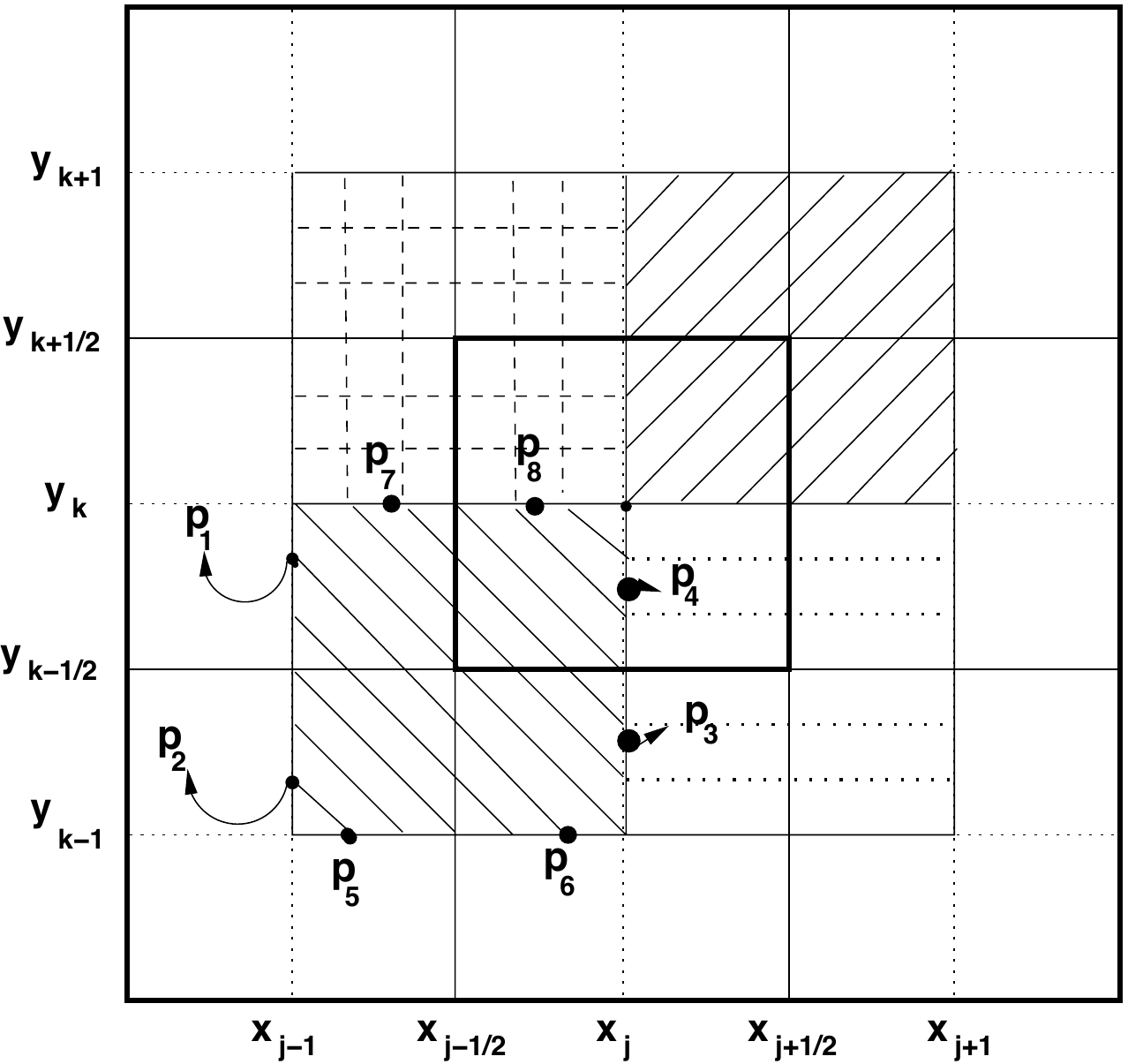}}
\caption{\label{fig:nt2d} SD2-2D: the construction of the two-dimensional grid}
\end{figure}

$(3)$ The NT2D scheme computes four averaged solutions at time step $\tmk + \Delta t_c$
(See the hachured regions in Figure \ref{fig:nt2d} corresponding to the cells
$C_{j-1/2,k-1/2}$,
$C_{j+1/2,k-1/2}$, $C_{j-1/2,k+1/2}$ and $C_{j+1/2,k+1/2}$
on the staggered grid). These averaged staggered solution  
are then properly projected back onto the original non-staggered grid
to obtain the desired solution (See \cite{NT2D}).

Repeating the same ideias presented by Rusanov in his modification of
Lax-Friedrichs' method, we compute the local speed of propragation at each
Riemann Problem. These local speeds define new non-uniform cells where
the evolution step will take place. 

{\bf Computing the local speed of propagation: } we begin with the cell   
$C_{j-1/2,k-1/2}$ to find the local speeds at the following Riemann Problems:

\begin{enumerate}
\item  \textbf{Y direction}:
\begin{enumerate} 
\item $$
\left\{
\begin{array}{ll}
\overline{S}^\kappa_{j-1,k-1}, & x_{j-3/2} \le x \le
x_{j-1/2}, \quad y_{k-3/2} \le y \le
y_{k-1/2} \\ 
\overline{S}^\kappa_{j-1,k}, & x_{j-3/2} \le x \le
x_{j-1/2}, \quad y_{k-1/2} \le y \le
y_{k+1/2}.   
\end{array}
\right.
$$
The local speed  $a^{y}_{j-1,k-1/2}$ defines the following points:
\begin{align*}
&p_1  =  (x_{j-1},\,y_{k-1/2} - a^y_{j-1,k-1/2}\dtcfl), \\
&p_2  = (x_{j-1},\,y_{k-1/2} + a^y_{j-1,k-1/2}\dtcfl)
\end{align*}
sketched in Figure \ref{fig:nt2d}. We denote the distance
between them by $\Delta y_{j-1,k-1/2} :=
2a^y_{j-1,k-1/2}\dtc$. 
\item  
$$
\left\{
\begin{array}{ll}
\overline{S}^\kappa_{j,k-1}, & x_{j-1/2} \le x \le
x_{j+1/2}, \quad y_{k-3/2} \le y \le
y_{k-1/2} \\ 
\overline{S}^\kappa_{j,k}, & x_{j-1/2} \le x \le
x_{j+1/2}, \quad y_{k-1/2} \le y \le
y_{k+1/2}.   
\end{array}
\right.
$$
The local speed
$a^{y}_{j,k-1/2}$ defines the points 
\begin{align*}
&p_3  = (x_{j},\,y_{k-1/2} - a^{y}_{j,k-1/2}\dtcfl), \\
&p_4  = (x_{j},\,y_{k-1/2} + a^{y}_{j,k-1/2}\dtcfl)
\end{align*}
shown in Figure \ref{fig:nt2d} and the distance between them is 
$\Delta y_{j,k-1/2} = 2a^{y}_{j,k-1/2}\dtcfl$.
\end{enumerate}
\item \textbf{X direction}:
\begin{enumerate}
\item
$$
\left\{
\begin{array}{ll}
\overline{S}^\kappa_{j-1,k-1}, & x_{j-3/2} \le x \le
x_{j-1/2}, \quad y_{k-3/2} \le y \le
y_{k-1/2} \\ 
\overline{S}^\kappa_{j,k-1}, & x_{j-1/2} \le x \le
x_{j+1/2}, \quad y_{k-3/2} \le y \le
y_{k-1/2}.   
\end{array}
\right.
$$
The local speed
$a^{x}_{j-1/2,k-1}$ defines the points
\begin{align*}
&p_5  = (x_{j-1/2} - a^{x}_{j-1/2,k-1}\dtcfl,\,y_{k-1}),\\
&p_6  = (x_{j-1/2} + a^{x}_{j-1/2,k-1}\dtcfl,\,y_{k-1})
\end{align*}
also shown in \ref{fig:nt2d} and $\Delta
x_{j-1/2,k-1}
= 2a^{x}_{j-1/2,k-1}\dtcfl$ is the distance between them.
\item 
$$
\left\{
\begin{array}{ll}
\overline{S}^\kappa_{j-1,k}, & x_{j-3/2} \le x \le
x_{j-1/2}, \quad y_{k-1/2} \le y \le
y_{k+1/2} \\ 
\overline{S}^\kappa_{j,k}, & x_{j-1/2} \le x \le
x_{j+1/2}, \quad y_{k-1/2} \le y \le
y_{k+1/2}.   
\end{array}
\right.
$$
The local speed
$a^{x}_{j-1/2,k}$ defines the points
\begin{align*}
p_7  = (x_{j-1/2} - a^{x}_{j-1/2,k}\dtcfl,\,y_{k}), \\
p_8  = (x_{j-1/2} + a^{x}_{j-1/2,k}\dtcfl,y_{k})
\end{align*}                         
shown in Figure \ref{fig:nt2d} and
$\Delta X_{j-1/2,k} = 2a^{x}_{j-1/2,k}\dtcfl$ denotes
the distance between them.
\end{enumerate}
\end{enumerate}

\begin{figure}[h]
\centerline{\includegraphics[scale=0.7]{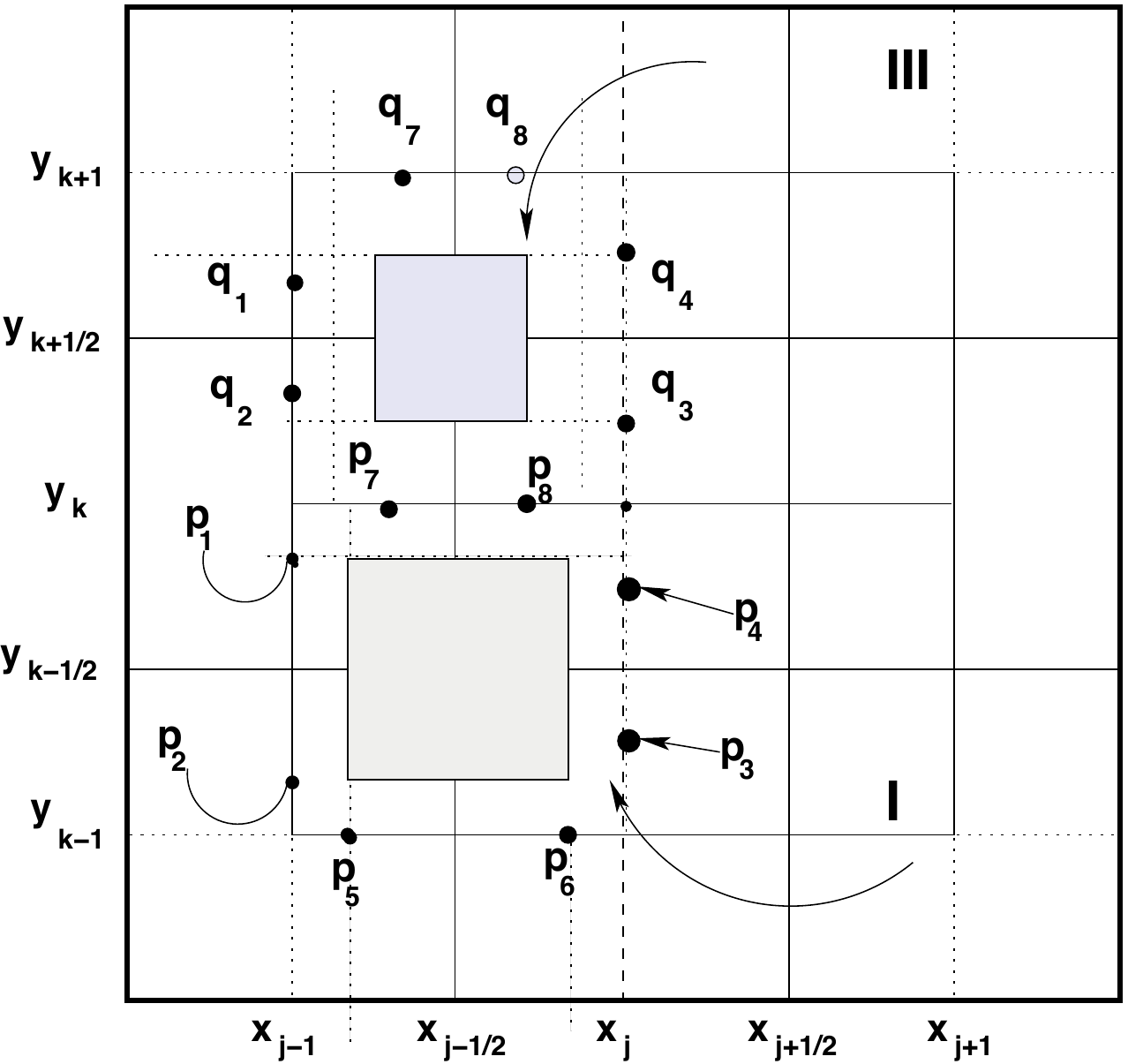}}
\caption{\label{fig:malhaktbi2}Regions I and III}
\end{figure}

Given these four local speed of wave propagation $a^y_{j-1,k-1/2}$,
$a^y_{j,k-1/2}$, $a^{x}_{j-1/2,k-1}$ e $a^{x}_{j-1/2,k}$, we can define
the Region I (also represented by $R_{j-1/2,k-1/2}$) as follows:  
\begin{itemize}
\item[$\rhd$] \textbf{Region I:}
  \begin{eqnarray*}
    \label{regionI}
    R_{j-1/2,k-1/2} &:=& [x_{j-1/2} -
    b^x_{j-1/2,k-1/2}\dtcfl, \, x_{j-1/2} +
    b^x_{j-1/2,k-1/2}\dtcfl] \times 
    \\
    & & [y_{k-1/2} -
    b^y_{j-1/2,k-1/2}\dtcfl, \, y_{k-1/2} +
    b^y_{j-1/2,k-1/2}\dtcfl] \\ \\ 
    \mbox{where} & &\hspace{-0.5cm}  
    b^x_{j-1/2,k-1/2}   :=  \max\{a^x_{j-1/2,k},
    a^x_{j-1/2,k-1}\} \hspace{0.3cm}\mbox{e}\hspace{0.3cm} \\ 
    & & \hspace{-0.5cm} b^y_{j-1/2,k-1/2}  :=  \max\{a^y_{j,k-1/2},
    a^y_{j-1,k-1/2}\}.
  \end{eqnarray*} 
\end{itemize}
 Figure \ref{fig:malhaktbi2} shows the new cell
$R_{j-1/2,k-1/2}$ of the new staggered non-uniform grid.

\begin{figure}[h]
\label{fig:malha_ktbi3}
\centerline{\includegraphics[scale=0.75]{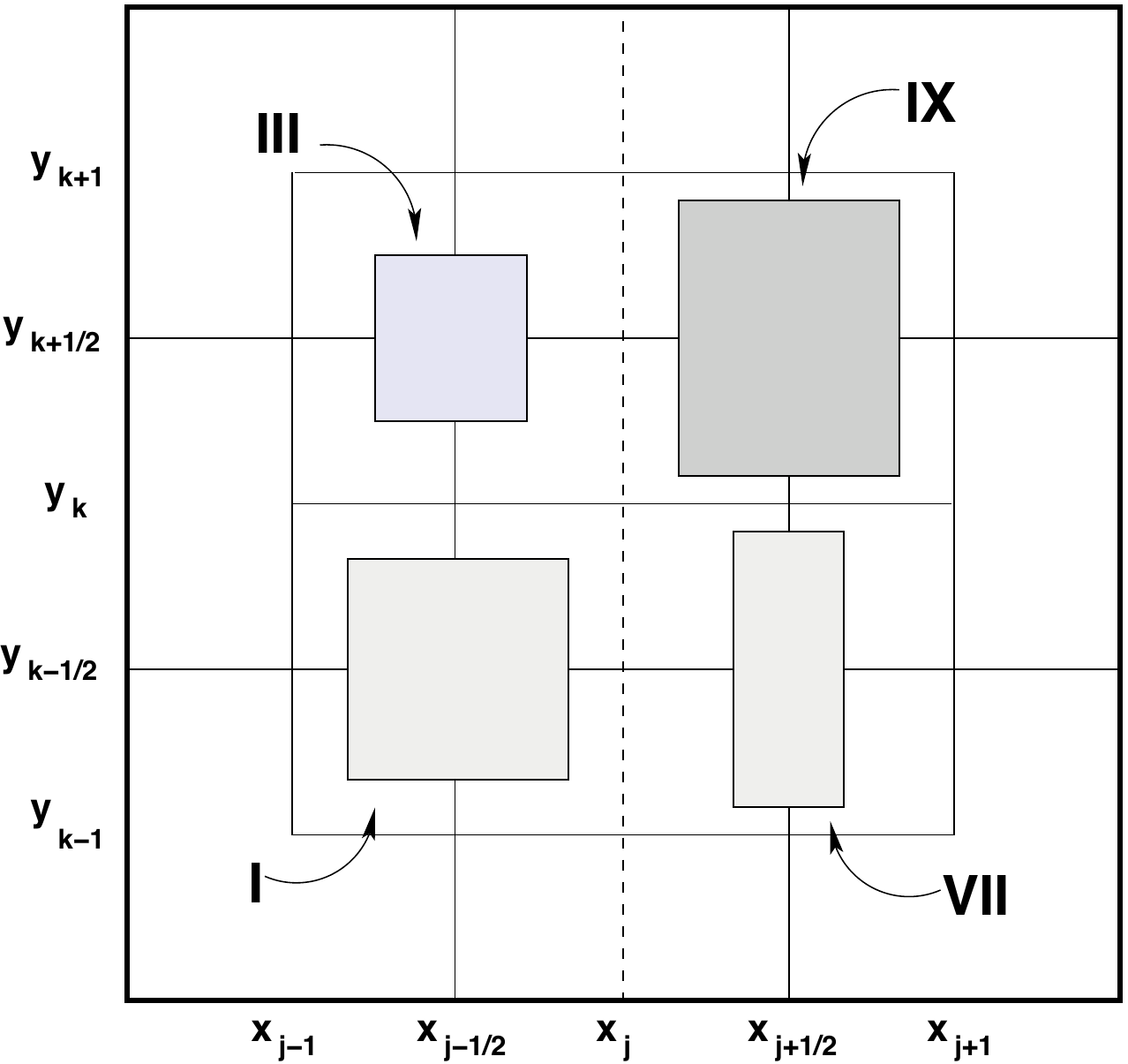}}
\caption{Regions I, III, VII and IX}
\end{figure}

We repeat the same procedures with the cell
$C_{j-1/2,k+1/2}$ to define the Region III in terms of the local speed of propagation $a^y_{j-1,k+1/2}$,
$a^y_{j,k+1/2}$, $a^{x}_{j-1/2,k+1}$ e $a^{x}_{j-1/2,k}$. 
These local speeds determine analogously the points
$q_1$ a $q_8$ sketched in Figure
\ref{fig:malhaktbi2} and define the new Region III:

\begin{itemize}
\item[$\rhd$] \textbf{Region III:}
  \begin{eqnarray*}
    \label{regionIII}
    R_{j-1/2,k+1/2} & := & [x_{j-1/2} -
    b^x_{j-1/2,k+1/2}\dtcfl, \, x_{j-1/2} +
    b^x_{j-1/2,k+1/2}\dtcfl] \times \nonumber \\
    & & [y_{k+1/2} -
    b^y_{j-1/2,k+1/2}\dtcfl, \, y_{k+1/2} +
    b^y_{j-1/2,k+1/2}\dtcfl] \\ \\ 
    \mbox{where} & &\hspace{-0.8cm}
    b^x_{j-1/2,k+1/2}  :=  \max\{ a^x_{j-1/2,k},
    a^x_{j-1/2,k+1}\} \hspace{0.3cm}\mbox{and}\hspace{0.3cm}\\  
    & &\hspace{-0.8cm} b^y_{j-1/2,k+1/2}  :=\max\{ a^y_{j,k+1/2},a^y_{j-1,k+1/2}\}.
  \end{eqnarray*}
\end{itemize}
Following these same procedures with the staggered cells $C_{j+1/2,k-1/2}$ and $C_{j+1/2,k+1/2}$ we define Regions VII and IX shown in Figure \ref{fig:malha_ktbi3}. 
We will denote by Group A the set of Regions I, III, VII and IX.
Finally, to finish the construction of the new non-uniform staggered grid, we need to define:
\begin{itemize}
\item Four cells which lie in the empty spaces between the regions of Group A. 
This new set will be denoted by Group B.
\item  A central region where the solution is smooth. 
\end{itemize}

The definitions of the cells of Group B will determine the central region. This central region will not be a rectangle, but a set of retangles. There are many ways to define the cells of
Group B. Our definition will be as follows (See Figura \ref{fig:malha_ktbi4}):
\begin{figure}[h]
\centerline{\includegraphics[scale=0.75]{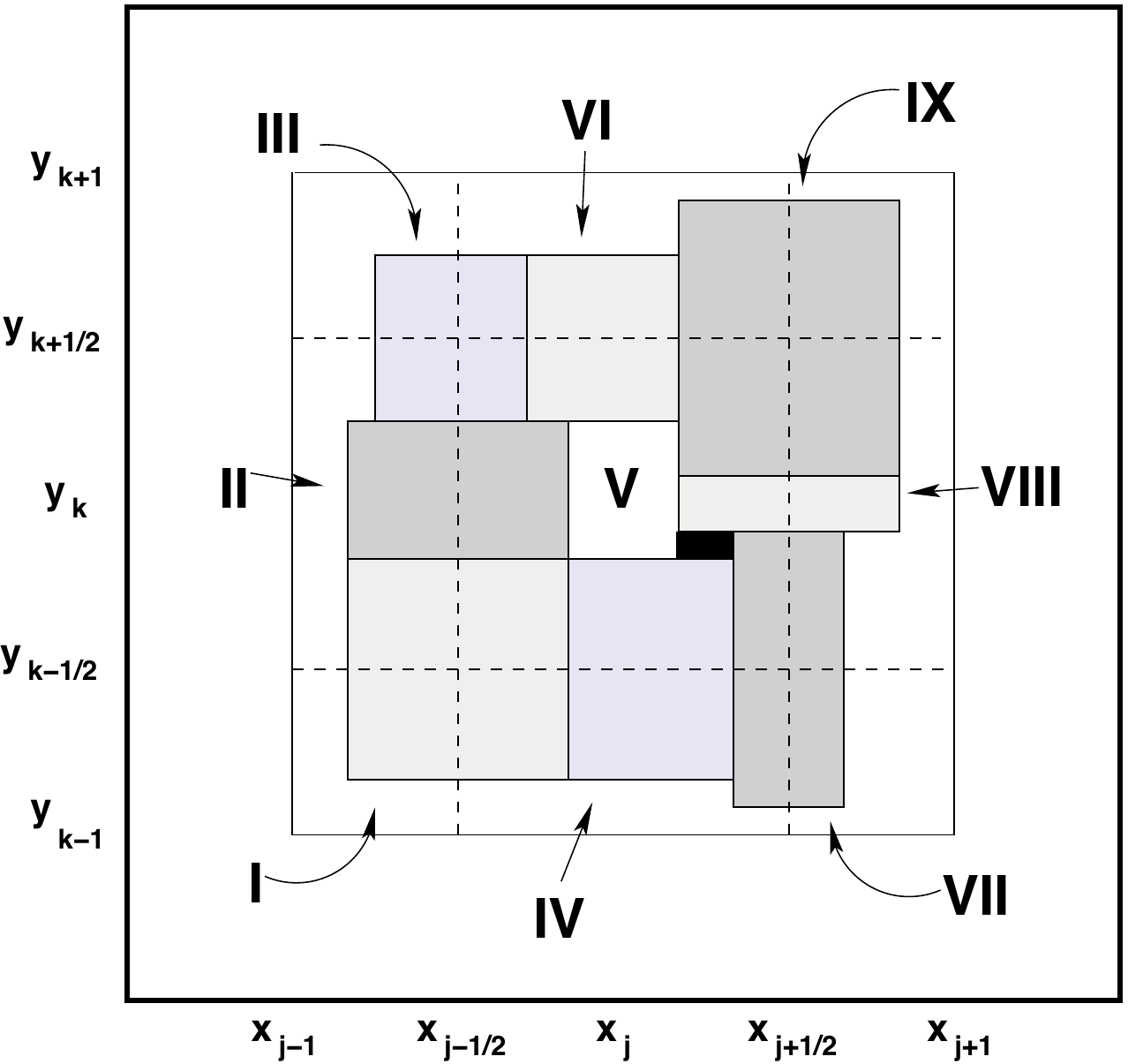}}
\caption{SD2-2D: The two dimensional semi-discrete central scheme SD2-2D: construction of the non-uniform staggered grid.}
\label{fig:malha_ktbi4}
\end{figure}
\begin{itemize}
\item[$\rhd$] \textbf{Region II:}
  \begin{eqnarray*}
    \label{regionLaterais}
    R_{j-1/2,k} &:=& [x_{j-1/2} -
    c^x_{j-1/2,k}\dtcfl, \, x_{j-1/2} +
    c^x_{j-1/2,k}\dtcfl] \times \\
    & & [y_{k-1/2} +
    b^y_{j-1/2,k-1/2}\dtcfl, \, y_{k+1/2} -
    b^y_{j-1/2,k+1/2}\dtcfl]  \\ \\ 
    & & \mbox{where}\quad
    c^x_{j-1/2,k}   :=  \max\{b^x_{j-1/2,k-1/2},
    b^x_{j-1/2,k+1/2}\}
  \end{eqnarray*}
\end{itemize}
The Regions VI and VIII can be obtained analogously 
As soon as the cells of Group A and B are determined, the central
region is automatically defined:
\begin{itemize}
\item[$\rhd$] \textbf{Region V:}
  \begin{eqnarray*}
    R_{j,k}& := &[x_{j-1/2} +
    c^x_{j-1/2,k}\dtcfl, \, x_{j+1/2} -
    c^x_{j+1/2,k}\dtcfl] \times 
     \\
     & &
    [y_{k-1/2} +
    d^y_{j,k-1/2}\dtcfl, \, y_{k+1/2} -
    d^y_{j,k+1/2}\dtcfl]. 
  \end{eqnarray*}
\end{itemize}

We also would like to emphasize that our choice for these regions
does not introduce
more numerical diffusion. 
We will call BR the black rectangle that can be seen in
Figure \ref{fig:malha_ktbi4} and we notice that, by construction, its area
is proportional to $(\Delta t_c)^2$.

\subsection{The new SD2-2D central scheme using Algorithm  REA}
\label{sec:KT2d_REA}

After defining the new control volumes performed in the section above, we are
now able to develop our new SD2-2D central scheme  following the REA algorithm. 

\textbf{Reconstruction step: } 
We suppose that we know an approximated solution constant in each cell at time step $t^\kappa$ 
as in Equation  \eqref{eq:cellaver}.   
This approximated solution is then reconstructed as a piecewise
bilinear polinomial as defined in Equations 
 \eqref{eq:bilinear_reconstruction} and
\eqref{eq:dev_2ordem_2D}.
\textbf{Evolution step: } Let $D$ represents one of the nine regions defined above
 $R_{j\pm 1/2, k\pm 1/2}$, $R_{j, k\pm 1/2}$, $R_{j\pm 1/2, k}$, the
central region $R_{j,k}$ or the black rectangle BR. We will denote by $D^+$ the part of region 
$D$ inside the non-staggered cell $I_{j,k}$ and by
$D^-$, the part of region $D$ outside the cell
$I_{j,k}$. 

We integrate the conservation law  \eqref{eq:sateq} in the control volumes
$D \times [\tmk, \tmk + \dtcfl]$ to obtain an approximate averaged solution
  $\bar{w}^{\kappa +1}(D)$ at the next time step,
in each cell $D$ of the staggered non-uniform grid.

\textbf{Projection step: } These averaged solutions
$\overline{w}^{\kappa+1}(D)$ are then reconstructed as piecewise
bilinear polynomials
$\widetilde{w}^{\kappa+1}(x,y)$ in each of the ten regions 
$D$. These new reconstructions are then projected back onto
the original grid of uniform non-staggered cells,
\begin{equation}
  \label{eq:media_futuro}
\overline{S}^{\kappa+1}_{j,k} := \frac{1}{\Delta x \,\Delta y}
  \int_{\bigcup D} \widetilde{w}^{\kappa+1} (x,y) \,dx\,dy.
\end{equation}

The new reconstructions $\widetilde{w}^{\kappa+1} (x,y)$ are
defined analogously as in Equation~\eqref{eq:bilinear_reconstruction}. 
For instance, for Region $R_{j+1/2,k}$,
\begin{align}
  \label{eq:linear-approx-2d-future_jmaisk}
  {\widetilde{w}^{\kappa+1}_{j+1/2,k}}(x,y) = \overline{w}^{\kappa+1}_{j+1/2,k} +
  (w_x)^{\kappa+1}_{j+1/2,k}(x-x_{j+1/2}) + &
  (w_y)^{\kappa+1}_{j+1/2,k}(y-y_k),  
  \\ & (x,y) \in R_{j+1/2,k}. \nonumber
\end{align}
The numerical derivatives $(w_x)^{\kappa+1}_{j+1/2,k}$ and
$(w_y)^{\kappa+1}_{j+1/2,k}$ satisfy the conditions 
\begin{eqnarray}
  \label{eq:dev_2ordem_2D_x}
  (w_x)^{\kappa+1}_{j+1/2,k} & = & \left.\frac{\partial w}{\partial
    x}\,\right|_{(x_{j+1/2},y_k,t^{\kappa+1})} +
\mathcal{O}(\Delta x);  \\
  \label{eq:dev_2ordem_2D_y}
(w_y)^{\kappa+1}_{j+1/2,k} & = & \left.\frac{\partial w}{\partial
    y}\,\right|_{(x_{j+1/2},y_k,t^{\kappa+1})} + \mathcal{O}(\Delta y);  
\end{eqnarray}
in order to guarantee the second order approximation.
Also, the reconstruction $\widetilde{w}^{\kappa+1}_{j+1/2,k}(x,y)$
retains the conservation property 
\eqref{eq:tilde_preserva_conserv}.  We remark that this is a theoretical step
and it will not be necessary to compute these numerical derivatives 
in the final semi-discrete formulation.

This completes the construction of our totally discrete central scheme in a 
rectangular grid.  It is very laborious to write a totally discrete version of this
central scheme. Instead, we will proceed directly to our semi-discrete formulation.
In order to to this, we compute the following limit when $\dtcfl \rightarrow 0$,
\begin{align}
  \label{eq:semi-discrete_2D}
   \lim_{\dtcfl \ra 0} &
  \frac{\overline{S}_{j,k}(t+\dtcfl) - \overline{S}_{j,k}(t)}{\dtcfl} 
= \frac{d}{dt}\overline{S}_{j,k}(t) = 
  \nonumber \\
& = \lim_{\dtcfl \ra 0} \frac{1}{\dtc}\cdot
\frac{1}{\Delta x\, \Delta y}\left\{ 
\sum_{p = j\pm 1/2}\int_{R^+_{p,k+1/2}} \hspace{-0.6cm}
  \widetilde{w}^{\kappa+1}_{p,k+1/2} (x,y) \right. \nonumber \\
& \left. +
\sum_{p = j\pm 1/2}\int_{R^+_{p,k-1/2}} \hspace{-0.6cm}
  \widetilde{w}^{\kappa+1}_{p,k-1/2} (x,y)
+  \sum_{p = j\pm 1/2}\int_{R^+_{p,k}} \hspace{-0.2cm}
  \widetilde{w}^{\kappa+1}_{p,k} (x,y)\,
  dx\,dy \right. \nonumber \\
& \left.
  + \sum_{q = k\pm 1/2} \int_{R^+_{j,q}} \hspace{-0.4cm}
  \widetilde{w}^{\kappa+1}_{j,q} (x,y)\,
  dx\,dy   
 + \int_{R_{j,k}} \hspace{-0.2cm}\widetilde{w}^{\kappa+1}_{j,k} (x,y)
  \,dx\,dy \right.\nonumber \\  
 & \left. + \int_{RP}
  \widetilde{w}^{\kappa+1}_{j+1/2,k-1/2} (x,y) dxdy
- (\Delta x\, \Delta y)\overline{S}_{j,k}(t) \right\} 
\end{align}

The conservation property of the reconstructions 
$\widetilde{w}^{\kappa+1}_{j,k}$ in he regions $D$
results 
\begin{equation}
  \label{eq:cons-property}
  \int_{D} \widetilde{w}^{\kappa+1}(x,y)\, dx\,dy = |D|\cdot\overline{w}^{\kappa+1}(D),
\end{equation}
where $\overline{w}^{\kappa+1}(D)$ is the averaged solution in region 
$D$. Note that, by reconstruction, the area of regions
 I, III, VII and IX are proportional to $(\dtcfl)^2$, that is,
$$
\begin{array}{ll}
 \mbox{Region I: } &  |R_{j-1/2,k-1/2}| = \mathcal{O}\big((\dtcfl)^2\big) \\
 \mbox{Region III: } &  |R_{j- 1/2,k+ 1/2}| = \mathcal{O}\big((\dtcfl)^2\big) \\
 \mbox{Region VII: } &  |R_{j+ 1/2,k- 1/2}| =
\mathcal{O}\big((\dtcfl)^2\big) \\
 \mbox{Region IX: } &  |R_{j+ 1/2,k+ 1/2}| =
\mathcal{O}\big((\dtcfl)^2\big) \\
 \mbox{Region RP: } &  |RP| =
\mathcal{O}\big((\dtcfl)^2\big) \\
\end{array}
$$
For example, considering Region IX, we conclude
\begin{align}
 \label{eq:area_nula}
 & \int_{R_{j+1/2, k+1/2}} \hspace{-0.5cm}
  \widetilde{w}^{\kappa+1}_{j+1/2,k+1/2} (x,y)\, dx\,dy  =
  \mathcal{O}(\dtcfl)^2 \nonumber \\
  &\hspace{3.0cm} \Rightarrow  \lim_{\dtcfl \ra 0} \frac{1}{\dtcfl}
  \int_{R_{j+1/2 , k+1/2 }} \hspace{-0.5cm}
  \widetilde{w}^{\kappa+1}_{j+1/2 ,k+1/2} (x,y)\, dx\,dy = 0;
 \end{align}

Our goal is to obtain the semi-discrete formulation of this
new central cheme. Therefore, the Equation \eqref{eq:area_nula}
show that we do not need to computed the averaged soltions over
the Regions I, III, VII and IX. Note that, by simetry, we only need
to compute the solutions over the Regions VI, VIII, V and the Black Rectangle. For the 
Region   $R_{j+1/2,k}$, we obtain:
 \begin{align}
\label{regiao_X}
 & \int_{R^+_{j+1/2 ,k}} 
 \widetilde{w}^{\kappa+1}_{j+1/2,k} (x,y)
  dxdy  =  \nonumber\\ 
& =   \int_{R^+_{j+1/2,k}} \hspace{-0.0cm}
\Big[ \bigg(\overline{w}^{\kappa+1}_{j+1/2,k} +
  (w_x)^{\kappa+1}_{j+1/2,k}\cdot (x-x_{j+1/2}) +
  (w_y)^{\kappa+1}_{j+1/2,k}\cdot(y-y_k)\bigg)\Big] dxdy \nonumber \\
 & =   
  \overline{w}^{\kappa+1}_{j+1/2,k} \cdot |R^+_{j+1/2,k}| +
  \mathcal{O}((\dtcfl)^2).
\end{align}
Analogously, we compute the averaged solution over Region $R^+_{j,k+1/2}$:
\begin{align}
\label{regiao-Y}
 & \int_{R^+_{j,k+1/2 }} \hspace{-0.3cm}
 \widetilde{w}^{\kappa+1}_{j,k+1/2} (x,y)
  dxdy  =   
  \overline{w}^{\kappa+1}_{j,k+1/2} \cdot |R^+_{j,k+1/2 }| +
  \mathcal{O}((\dtc)^2).
\end{align}
Note that the solution has no discontinuities inside  $R_{j,k}$.  
So, it isn't necessary to reconstruction as a piecewise bilinear polynomials. 
The averaged solution is
\begin{align}
\label{regiao-inside}
 & \int_{R_{j,k}} \hspace{-0.2cm}
 \widetilde{w}^{\kappa+1}_{j,k} (x,y)
  dxdy  =  
  \overline{w}^{\kappa+1}_{j,k} \cdot |R_{j,k}|.
\end{align}

Substituting the Equations~\eqref{eq:area_nula}, \eqref{regiao_X},
 \eqref{regiao-Y} and \eqref{regiao-inside} in Equation
\eqref{eq:semi-discrete_2D}, we obtain:
\begin{align}
  \label{eq:semi-discrete2}
  \frac{d}{dt}\overline{S}_{j,k}(t) =& \lim_{\dtcfl \ra 0}
  \left\{\frac{c^x_{j-1/2,k}}{\Delta x}
    \overline{w}_{j-1/2,k}(t + \dtcfl) +
    \frac{c^x_{j+1/2,k}}{\Delta x} 
    \overline{w}_{j+1/2,k}(t + \dtcfl) \right.
    \nonumber \\
    & \left. +\frac{d^y_{j,k-1/2}}{\Delta y}
     \overline{w}_{j,k-1/2}(t + \dtcfl) +
    \frac{d^y_{j,k+1/2}}{\Delta Y} 
    \overline{w}_{j,k+1/2}(t + \dtcfl) \right.\nonumber\\
    &\left. - \left( \frac{d^y_{j,k+1/2} + d^y_{j,k-1/2}}{\Delta y}
      + \frac{c^x_{j+1/2,k} + c^x_{j-1/2,k}}{\Delta x}
    \right)\cdot \overline{w}_{j,k}(t + \dtcfl) \right.\nonumber \\
& \left.
+ \left( \frac{1}{\dtcfl}
\overline{w}_{j,k}(t + \dtcfl) -  \frac{1}{\dtcfl}
\overline{S}_{j,k}(t + \dtcfl) \right)\right\}. 
\end{align}

For the final formulation, we have to compute the averaged solution
over the non-uniform staggered grid.
\begin{equation*}
  \overline{w}^{\kappa+1}_{j+1/2,k},\overline{w}^{\kappa+1}_{j,k+1/2} \mbox{ e }\overline{w}^{\kappa+1}_{j,k},
\end{equation*}
To this end, we integrate the conservation law \eqref{eq:sist2_cont} over the control volumes 
\begin{equation*}
R_{j+1/2,k} \times [\tmk, \tmk+\dtcfl],\quad 
R_{j,k+1/2} \times [\tmk, \tmk+\dtcfl] \quad\mbox{e}\quad  
R_{j,k} \times [\tmk,\tmk+\dtcfl],
\end{equation*}
respectively.
Therefore, 
\begin{align}
\label{wjmaisk}
w^{\kappa+1}_{j+1/2,k} &   \nonumber 
 = \frac{1}{|R_{j+1/2,k}|}
\int_{R_{j+1/2,k}} s(x,y,t^{\kappa+1})\,dx\,dy \nonumber
\\
& =
\frac{1}{|R_{j+1/2,k}|}\int_{R_{j+1/2,k}} \widetilde{S}^\kappa(x,y)\, dx \,dy \nonumber \\
& -
\frac{1}{|R_{j+1/2,k}|}\int_{R_{j+1/2,k}}
\int_{t^\kappa}^{t^{\kappa+1}}
\bigg[\frac{\partial}{\partial
  x}\bigg(\vx(x,y,\tau)\cdot f(s(x,y,\tau))\bigg) \\
&\hspace{3.0cm}+ \frac{\partial}{\partial y}\bigg(
\vy(x,y,\tau)\cdot f(s(x,y,\tau)) \bigg) \bigg]\, dx\, dy
d\tau. \nonumber
\end{align}
Let us denote the  double integral 
 by Int$_{\widetilde{S}}$ and the flux integral by
 Int$_f$. The integral Int$_{\widetilde{S}}$ is computed analytically.
\begin{align}
\label{integral_espacial_1}
\mbox{Int}_{\widetilde{S}} = &
\frac{1}{2}(\overline{S}^\kappa_{j,k} +
\overline{S}^\kappa_{j+1,k}) \nonumber \\
& +\frac{1}{4}\left[\Delta x - c^x_{j+1/2,k}\dtcfl\right]
\cdot \left[ (S_x)^\kappa_{j,k} -
  (S_x)^\kappa_{j+1,k}\right] \\
& + \frac{\dtcfl}{4}\left[b^y_{j+1/2,k-1/2}-
  b^y_{j+1/2,k+1/2}\right] \cdot \left[ (S_y)^\kappa_{j,k} -
  (S_y)^\kappa_{j+1,k}\right]. \nonumber
\end{align}
To compute the flux integral, Int$_f$, we first denote the limits
of region $R_{j+1/2,k}$ as follows:
\begin{align*}
& a:= x_{j+1/2} - c^x_{j+1/2,k}\dtcfl \\
& b:= x_{j+1/2} + c^x_{j+1/2,k}\dtcfl \\
& c:= y_{k+1/2} + b^y_{j+1/2,k-1/2}\dtcfl \\
& d:= y_{k+1/2} - b^y_{j+1/2,k+1/2}\dtcfl \\
\end{align*}
Using the Calculus Fundamental Theorem together with the trapezoid rule, we obtain 
\begin{eqnarray}
\label{trapezio}
\mbox{Int}_f & = &  \frac{1}{2\Delta x_{j+1/2,k}}\int_{t^\kappa}^{t^{\kappa+1}}
\Big[\vx(b,d,\tau)\,f(s(b,d,\tau)) 
- \vx(a,d,\tau)\, f(s(a,d,\tau))
\nonumber \\
& & 
+ \vx(b,c,\tau)\, f(s(b,c,\tau)) 
- \vx(a,c,\tau)\, f(s(a,c,\tau))
\Big]\,d\tau 
\nonumber \\
& & + \frac{1}{2\Delta y_{j+1/2,k}} \int_{t^\kappa}^{t^{\kappa+1}} \Big[
\vy(b,d,\tau)\, f(s(b,d,\tau)) 
 - \vy(b,c,\tau)\, f(s(b,c,\tau))
\nonumber
\\
& & + \vy(a,d,\tau)\, f(s(a,d,\tau)) 
- \vy(a,c,\tau)\, f(s(a,c,\tau))\Big]\,d\tau 
\end{eqnarray}
If the CFL condition
\begin{equation}
  \label{eq:cflKT2D}
  \max\left(\frac{\dtcfl}{\Delta
      x}\max_S|\vx\, f'(s)|, \frac{\dtcfl}{\Delta
      y}\max_S|\vy\, f'(s)|\right)< \frac{1}{2}
\end{equation}
holds and since the functions $\vx\, f(s(x,y,\tau))$ and
$\vy\, f(s(x,y,\tau))$ are computed away from the discontinuities 
then they are differential functions of 
$\tau$ and therefore, the time integral can be approximated using the
middle point rule. Denoting  $t^{\kappa+1/2} :=t +
\dtcfl/2$, we obtain:
\begin{eqnarray}
\label{eq:ponto_meio}
\mbox{Int}_f &  = &  \displaystyle{\frac{1}{4c^x_{j+1/2,k}}
\Big[\vx(b,d,t^{\kappa + 1/2})\,f(S(b,d,t^{\kappa + 1/2})) 
 - \vx(a,d,t^{\kappa + 1/2})\, f(S(a,d,t^{\kappa + 1/2}))}
\nonumber \\
&  &
\displaystyle{+ \vx(b,c,t^{\kappa + 1/2})\, f(S(b,c,t^{\kappa + 1/2}))  - \vx(a,c,t^{\kappa + 1/2})\, f(S(a,c,t^{\kappa + 1/2}))
\Big]\,d\tau} 
\nonumber \\
& & \displaystyle{ + \frac{\alpha_Y}{2 - 2\alpha_Y(b^{y}_{j+1/2,k+1/2}
  -b^{y}_{j+1/2,k-1/2})}} \cdot \nonumber\\  
& & \displaystyle{\Big[\vy(b,d,t^{\kappa + 1/2})\, f(S(b,d,t^{\kappa + 1/2}))  
- \vy(b,c,t^{\kappa + 1/2})\, f(S(b,c,t^{\kappa + 1/2}))} \nonumber\\ 
& &\displaystyle{ + \vy(a,d,t^{\kappa + 1/2})\, f(S(a,d,t^{\kappa + 1/2})) 
 - \vy(a,c,t^{\kappa + 1/2})\, f(S(a,c,t^{\kappa + 1/2}))\Big]\,d\tau} 
\end{eqnarray}
where $\alpha_X =\dtcfl/\Delta X$ 
and $\alpha_Y = \dtcfl/\Delta Y$.

The midpoint  values are computed using the Taylor expansions and the
conservation law \eqref{eq:sateq}.
For instance,
\[
\left\{
\begin{array}{l}
\label{eq:midpoint_j+rku}
S(a,d,t^{\kappa + 1/2}) := S(a,d,t)   \nonumber
\\
 \hspace{3.0cm} \displaystyle{-\frac{\dtcfl}{2}
\Big( \vx(a,d,t)\, f(S(a,d,t)) \Big)_x} 
\displaystyle{- \frac{\dtcfl}{2}
\Big(\vy(a,d,t)\, f(S(a,d,t)) \Big)_y} 
\\
\\
\displaystyle{S(a,d,t)  := 
\overline{S}^\kappa_{j+1,k} - \Delta x\, (S_x)^\kappa_{j+1,k}
 \Big(\frac{1}{2} - \alpha_X c^x_{j+1/2,k}\Big)}  
\displaystyle{- \Delta y\, (S_y)^\kappa_{j+1,k}
\Big(\frac{1}{2} - \alpha_Y b^{y}_{j+1/2,k+1/2}\Big)}.
\end{array}
\right.
\]
As the time step $\dtc$ goes to zero, the limit
\begin{eqnarray}
\lim_{\dtcfl \ra 0} S(a,d,t^{\kappa + 1/2}) = \overline{S}^\kappa_{j+1,k} - \frac{\Delta X}{2}\, (S_x)^\kappa_{j+1,k} 
- \frac{\Delta Y}{2}\, (S_y)^\kappa_{j+1,k} := S^{++}_{j+1/2,k-1/2}.
\end{eqnarray}
These are called the intermediate values and their general form is
\begin{eqnarray}
 S^{\pm\pm}_{j+1/2,k+1/2} & = & \overline{S}^\kappa_{j+1/2\pm 1/2, k+1/2\pm 1/2} \pm \frac{\Delta X}{2}\, (S_x)^\kappa_{j+1/2\pm 1/2, k+1/2\pm 1/2}(x_{j+1/2} - x_{j+1/2\pm 1/2}) \nonumber \\
 & & \pm  \frac{\Delta Y}{2}\, (S_y)^\kappa_{j+1/2\pm 1/2, k+1/2\pm 1/2}(y_{k+1/2} - y_{k+1/2\pm 1/2})
\end{eqnarray}

We notice that the cell averages $w^{\kappa+1}_{j,k+1/2}$ and
$w^{\kappa+1}_{j,k}$ are obtained analogously to \eqref{wjmaisk}.
And also, $w^{\kappa+1}_{j-1/2,k} =
w^{\kappa+1}_{j+1/2-1,k}$ e $w^{\kappa+1}_{j,k-1/2} =w^{\kappa+1}_{j,k+1/2-1}$.

Substituting all these cell averages in time step $\tmk + \dtcfl$
into Equation \eqref{eq:semi-discrete2} and computing the limit when 
 $\dtcfl \ra 0$, we obtain the second order central scheme in
semi-discrete formulation:
\begin{equation}
\label{eq:semi_discreto_conservativo}
   \frac{d}{dt}\overline{S}_{jk}(t)  =
   -\frac{H^{x}_{j+1/2,k}(t)-H^{x}_{j-1/2,k}(t)}{\Delta X}
   -\frac{H^{y}_{j,k+1/2}(t)-H^{y}_{j,k-1/2}(t)}{\Delta Y},
\end{equation} 
where the numerical fluxes are 
\begin{subequations}\label{eq:fluxo_kt2D}
\begin{align}
H^{x}_{j+1/2,k}(t) = & \frac{1}{4}\Big\{ \vx_{j+1/2,k+1/2}(t)\Big[f(S^{+-}_{j+1/2,k+1/2}(t)) +
f(S^{--}_{j+1/2,k+1/2}(t))\Big] \nonumber \\
&  + \vx_{j+1/2,k-1/2}(t)\Big[f(S^{++}_{j+1/2,k-1/2}(t)) +
f(S^{-+}_{j+1/2,k-1/2}(t))\Big] \Big\} \nonumber \\
&  - \frac{c^x_{j+1/2,k}}{2}\Big[ S^{+}_{j+1/2,k}(t) -
S^{-}_{j+1/2,k}(t)\Big]; \label{eq:fluxo_x_2d} \\
H^{y}_{j,k+1/2}(t) & =  \frac{1}{4}\Big\{ \vy_{j+1/2,k+1/2}(t)\Big[f(S^{-+}_{j+1/2,k+1/2}(t)) +
f(S^{--}_{j+1/2,k+1/2}(t))\Big] \nonumber \\
& + \vy_{j-1/2,k+1/2}(t)\Big[f(S^{++}_{j-1/2,k+1/2}(t)) +
f(S^{+-}_{j-1/2,k+1/2}(t))\Big] \Big\}\nonumber \\
&  - \frac{d^y_{j,k+1/2}}{2}\Big[ S^{+}_{j,k+1/2}(t) - S^{-}_{j,k+1/2}(t)\Big].\label{eq:fluxo_y_2d}
\end{align}
\end{subequations}

If the numerical derivatives are equal to zero then we
obtain the two-dimensional Rusanov's central scheme in
semi-discrete formulation.
\begin{subequations}\label{eq:rusa_2d_fluxos}
\begin{align}
\mbox{{Rusa}}^{x}_{j+1/2,k}(t)  = & \frac{1}{4}\Big\{ \vx_{j+1/2,k+1/2}(t)\Big[f(S_{j+1,k}(t)) +
f(S_{j,k}(t))\Big] \nonumber \\
 & + \vx_{j+1/2,k-1/2}(t)\Big[f(S_{j+1,k}(t)) +
f(S_{j,k}(t))\Big] \Big\} \nonumber \\
 & - \frac{c^x_{j+1/2,k}}{2}\Big[ S_{j+1,k}(t) -
S_{j,k}(t)\Big]; \label{eq:fluxo_x_2d_rusanov} \\
\mbox{{Rusa}}^{y}_{j,k+1/2}(t)  = & \frac{1}{4}\Big\{ \vy_{j+1/2,k+1/2}(t)\Big[f(S_{j,k+1}(t)) +
f(S_{j,k}(t))\Big] \nonumber \\
 &+ \vy_{j-1/2,k+1/2}(t)\Big[f(S_{j,k+1}(t)) +
f(S_{j,k}(t))\Big] \Big\}\nonumber \\
 & - \frac{d^y_{j,k+1/2}}{2}\Big[ S_{j,k+1}(t) - S_{j,k}(t)\Big].\label{eq:fluxo_y_2d_rusanov}
\end{align}
\end{subequations}

This new two-dimensional semi-discrete central scheme with 
the numerical fluxes given by  \eqref{eq:fluxo_kt2D} 
or 
\eqref{eq:rusa_2d_fluxos} comprises a system of ordinary differential equations.
To solve this system, we use the explicit second order Runge-Kutta method.

\subsection{The velocity field } Finally, to
complete the description of the genuinely two-dimensional KT
scheme, we have to define the velocity field. The velocity
is defined at the vertices of the cells. We can not use
directly the velocity field from the Raviart-Thomas space as
we did in the dimension by dimension approach. Instead we
will use a bilinear interpolation of it preserving the null
divergence necessary for the incompressible flows. For
instance, to compute the value of $\vx_{j+1/2,k+1/2}$ on the
vertex $(x_{j+1/2},y_{k+1/2})$ at some time step $t^m$, we
have to use all the four cells which share this vertex,

\begin{eqnarray*}
  \vx_{j+1/2,k+1/2} & = & \frac{1}{2}\Big(\vx_{j,k+1/2} + \vx_{j+1,k+1/2}\Big)\\
  & = & \frac{1}{2}\Big[\frac{1}{2}\Big(\vx_{j,k} +
  \vx_{j,k+1}\Big) + \frac{1}{2}\Big(\vx_{j+1,k} +
  \vx_{j+1,k+1}\Big)\Big] \\ 
  & = & +
  \frac{1}{2}\Big[\frac{1}{2}\Big(\frac{1}{2}(\vec{v_r}_{j+1/2,k}
  - \vec{v_l}_{j-1/2,k}) \\ 
  & & + \frac{1}{2}(\vec{v_r}_{j+1/2,k+1} -
  \vec{v_l}_{j-1/2,k+1}\Big) + \frac{1}{2}  
  \Big(\frac{1}{2}(\vec{v_r}_{j+3/2,k} - \vec{v_l}_{j+1/2,k}) + \\
  & & \frac{1}{2}(\vec{v_r}_{j+3/2,k+1} - \vec{v_l}_{j+1/2,k+1})\Big) \\
  & = & \frac{1}{8}\Big(\vec{v_r}_{j+1/2,k} -
  \vec{v_l}_{j-1/2,k} + \vec{v_r}_{j+1/2,k+1} -  
  \vec{v_l}_{j-1/2,k+1} \\
  &  & + \vec{v_r}_{j+3/2,k} - \vec{v_l}_{j+1/2,k} +
  \vec{v_r}_{j+3/2,k+1} - \vec{v_l}_{j+1/2,k+1} 
  \Big)
\end{eqnarray*}

\section{TWO-DIMENSIONAL NUMERICAL EXPERIMENTS}
\label{sec:aplication}

We present and compare the results for numerical simulations
of two-dimensional, two-phase flow associated with two
distinct flooding problems using the KT scheme with the
dimension by dimension approach (KT dxd) and the genuinely
two-dimensional formulation (SD2-2D).  The first problem is
a two-dimensional flow in a rectangular heterogeneous
reservoir (called slab geometry) with $256$ m $\times$ $64$
m in size, and the second is a two-dimensional flow in a
5-spot geometry homogeneous reservoir having $64$ m $\times$
$64$ m.

In the 5-spot geometry homogeneous reservoir simulation we
used two distinct uniform five-spot well configurations
intended to illustrate different flow patterns, with
parallel and diagonal grid orientations, and boundary
behavior.

In all simulations the reservoir contains initially $79\%$
of oil and $21\%$ of water.  Water is injected at a constant
rate of $0.2$ pore volumes every year.
 
The fractional volumetric flow, the total mobility, and the
relative permeabilities are assumed to be:
\[
f(s) = \frac{k_{rw}(s)/\mu_w}{\lambda(s)},\quad
\lambda(s)=\frac{k_{rw}(s)}{\mu_w}+\frac{k_{ro}(s)}{\mu_o},
\]
and 
\[
k_{ro}(s) = (1 - (1-s_{ro})^{-1}s)^{2}, \quad
k_{rw}(s) = (1-s_{rw})^{-2}(s-s_{rw})^{2},
\]
where $s_{ro}=0.15$ and $s_{rw}=0.2$ are the residual oil
and water saturations, respectively. The viscosity of oil
and water used in all simulations are $\mu_o=10.0\,cP$ and
$\mu_o=0.05\,cP$.

For the heterogeneous reservoir studies we consider a scalar
absolute permeability field $K({\bf x})$ taken to be
log-normal (a fractal field, see \citep{GLPZ} and
\citep{fpcross} for more details) with moderately large
heterogeneity strength. The spatially variable permeability
field is defined on a 256 $\times$ 64 grid with three
different values of the coefficient of variation $C_v$
(standard deviation)/mean: $0.5, 1.0,$ and $2.4$.

\begin{figure}[htbp!]
\includegraphics[scale=0.3]{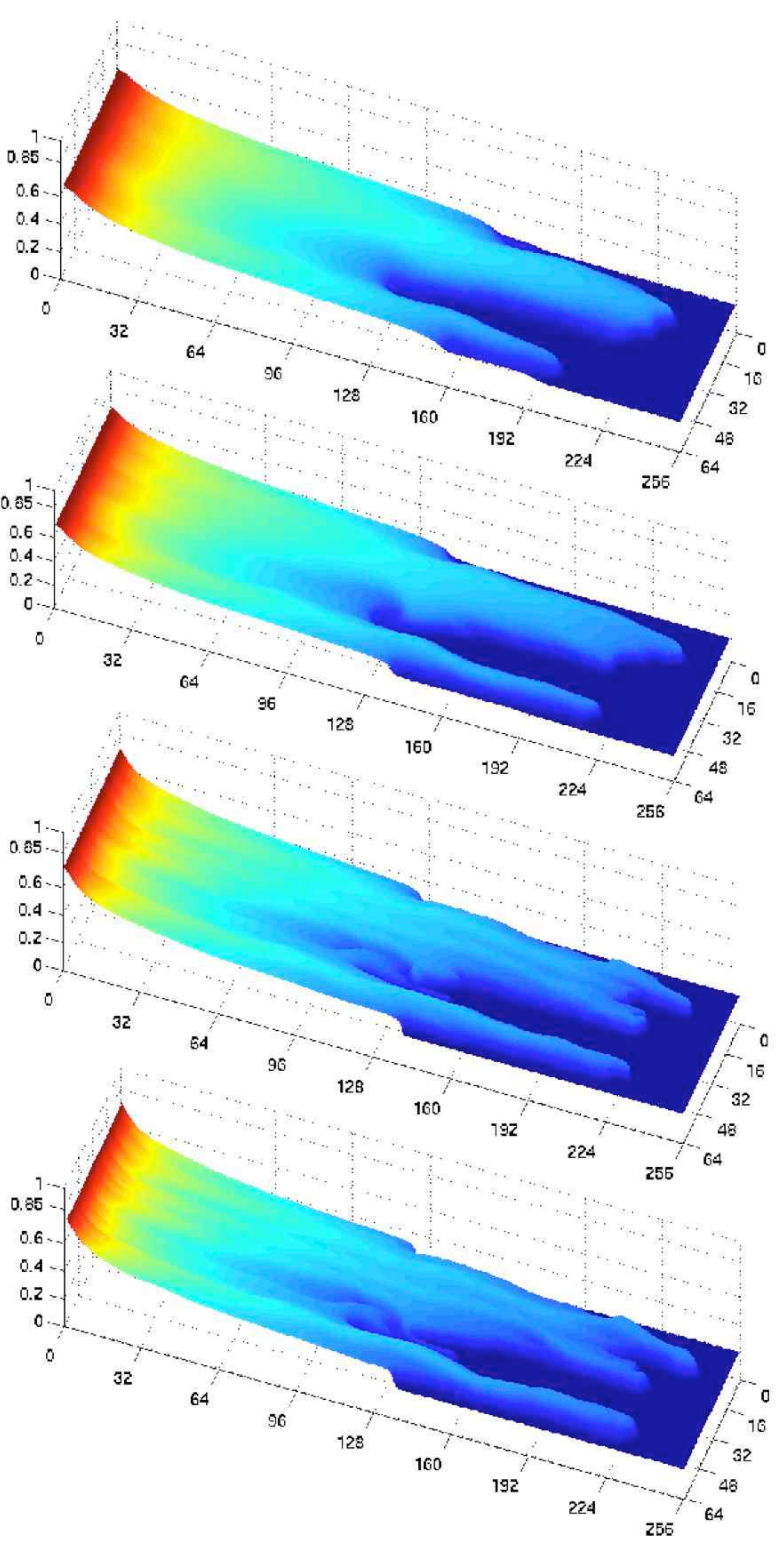}
\caption{Water saturation surface plots for two-phase flow
  in a two-dimensional heterogeneous reservoir having $256$
  m $\times$ $64$ m, with the coefficient of variation
  $C_v=0.5$ and viscous ratio 20. From top to bottom:$1)$
  (KT dxd) scheme with $256\times64$ grid; $2)$ (KT dxd)
  scheme with $512\times128$ grid; $3)$ (KT two) scheme with
  $256\times64$ grid.}
\label{slab}
\end{figure}

\begin{figure}[htbp!]
\includegraphics[scale=0.55]{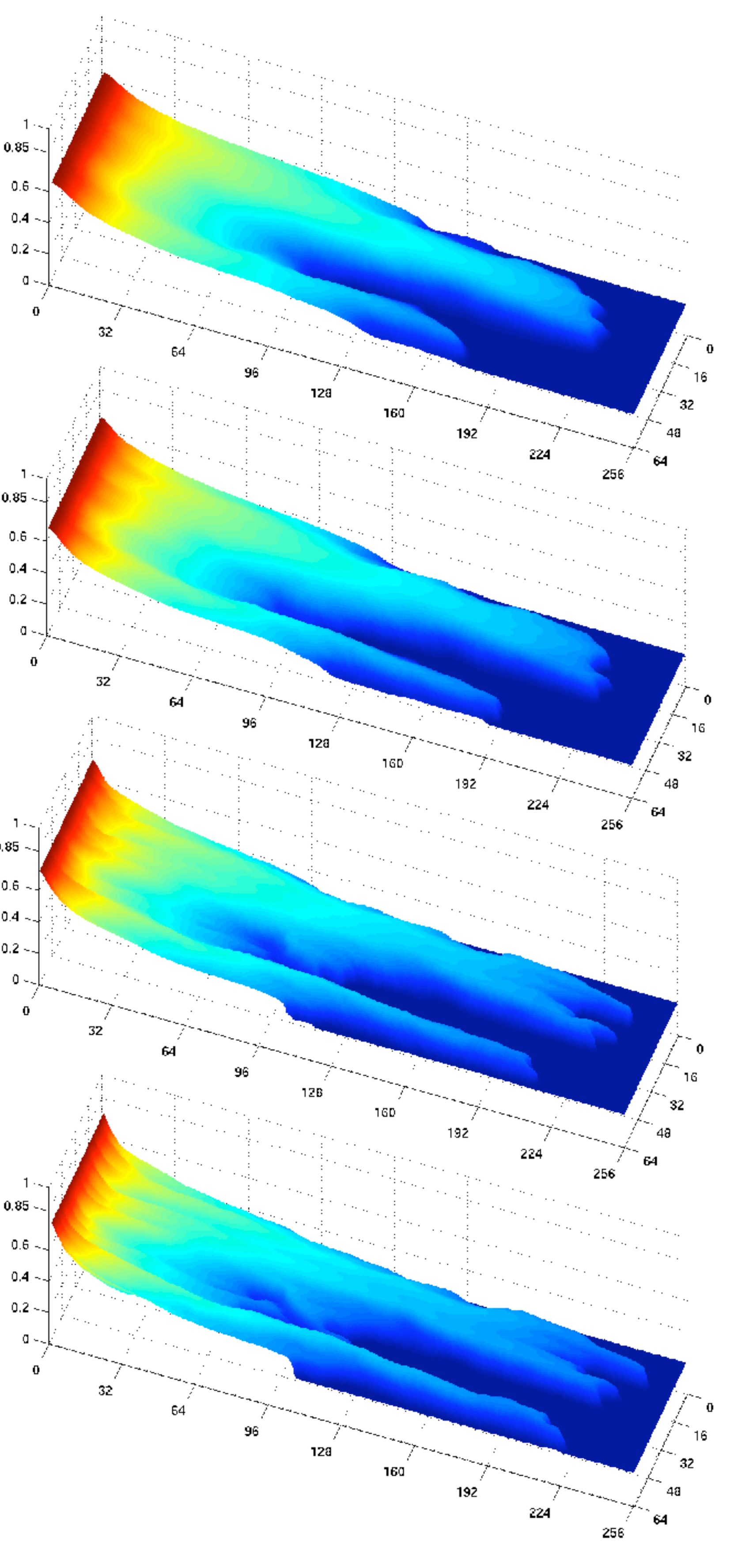}
\caption{Water saturation surface plots for two-phase flow
  in a two-dimensional heterogeneous reservoir having $256$
  m $\times$ $64$ m, with the coefficient of variation
  $CV=1.2$ and viscous ratio 20. From top to bottom:$1)$
  (KT dxd) scheme with $256\times64$ grid; $2)$ (KT dxd)
  scheme with $512\times128$ grid; $3)$ (KT two) scheme with
  $256\times64$ grid.}
\label{slab_CV1}
\end{figure}

\begin{figure}[htbp!]
\includegraphics[scale=0.55]{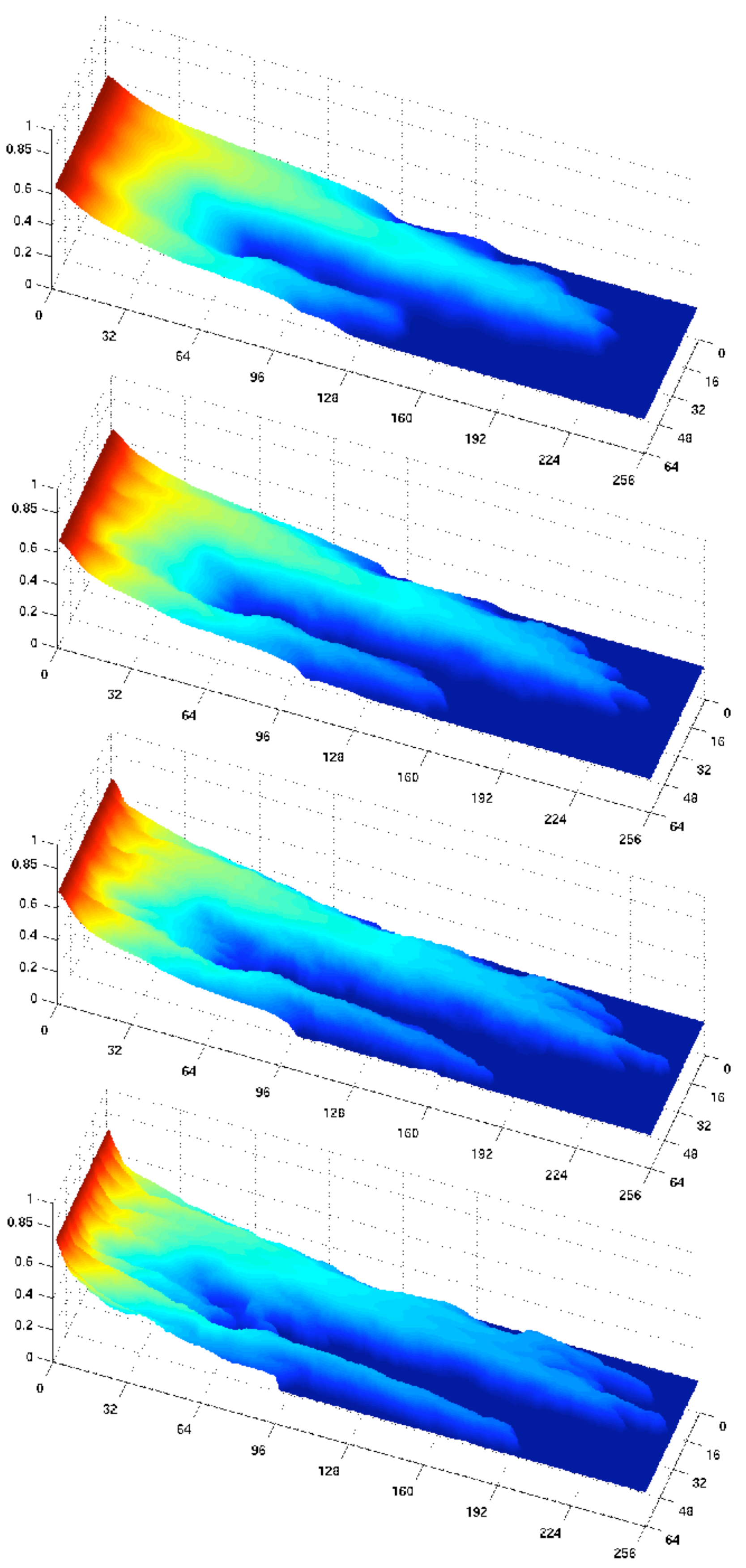}
\caption{Water saturation surface plots for two-phase flow
  in a two-dimensional heterogeneous reservoir having $256$
  m $\times$ $64$ m, with the coefficient of variation
  $CV=2.2$ and viscous ratio 20. From top to bottom:$1)$
  (KT dxd) scheme with $256\times64$ grid; $2)$ (KT dxd)
  scheme with $512\times128$ grid; $3)$ (KT two) scheme with
  $256\times64$ grid.}
\label{slab_CV2}
\end{figure}

The boundary conditions, injection and production
specifications for two-phase flow equations
\eqref{eq:preeq}-\eqref{eq:sateq}) are as follows. For the
horizontal slab geometry, injection is made uniformly along
the left edge of the reservoir and the (total) production
rate is taken to be uniform along the right edge; no flow is
allowed along the edges appearing at the top and bottom of
the reservoir in the graphics (Figures \ref{slab},
\ref{slab_CV1}, and \ref{slab_CV2}).

\begin{figure}[htbp]
\centering
\subfigure[\small NT2D: $64 \! \times \! 64$ células
]{\label{NT_dg}
\includegraphics[scale=0.35]{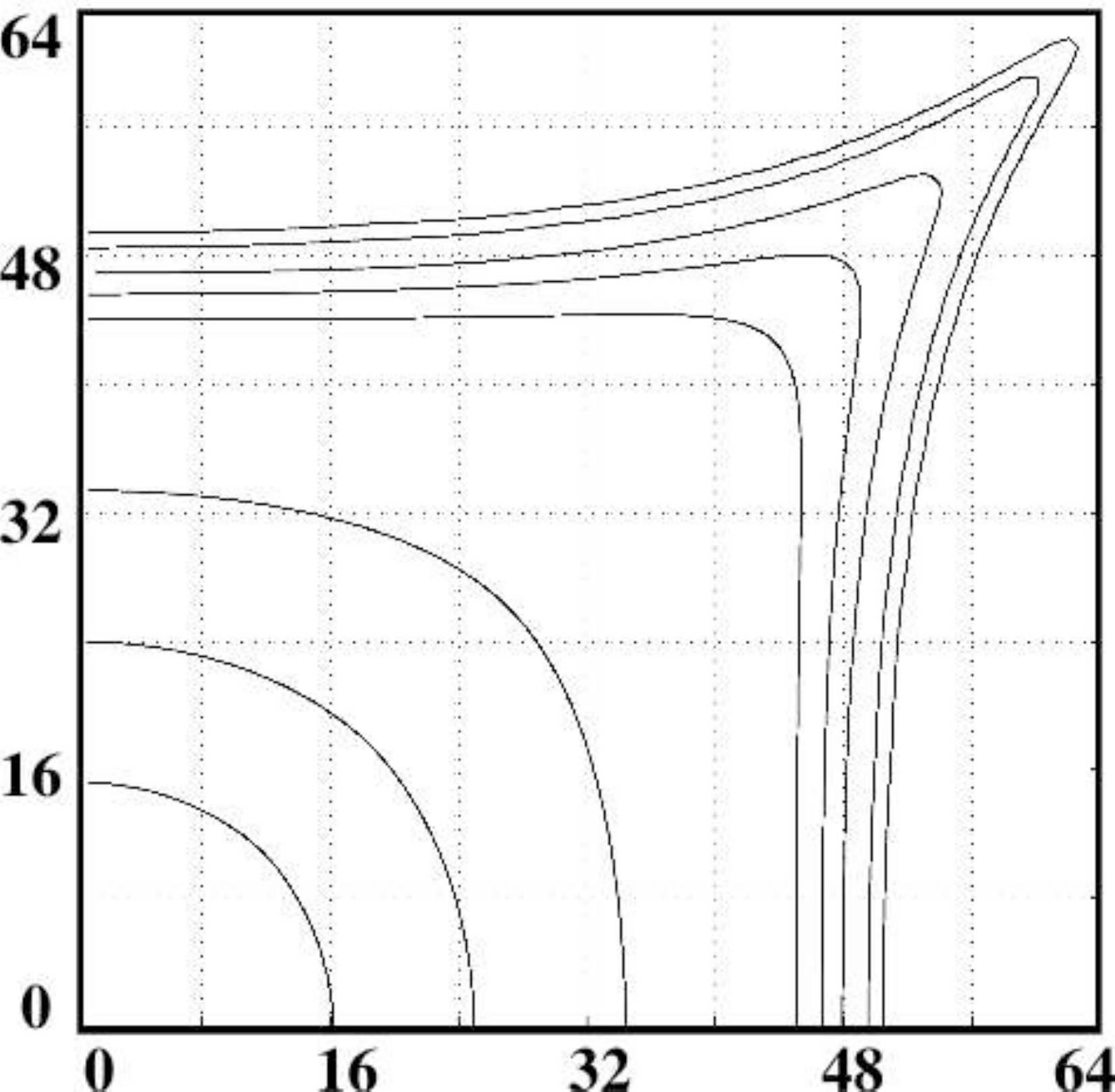}} 
\subfigure[\small KTdxd: $64 \! \times \! 64$ células
]{\label{KTdxd_pg}\quad\quad
\includegraphics[scale=0.7]{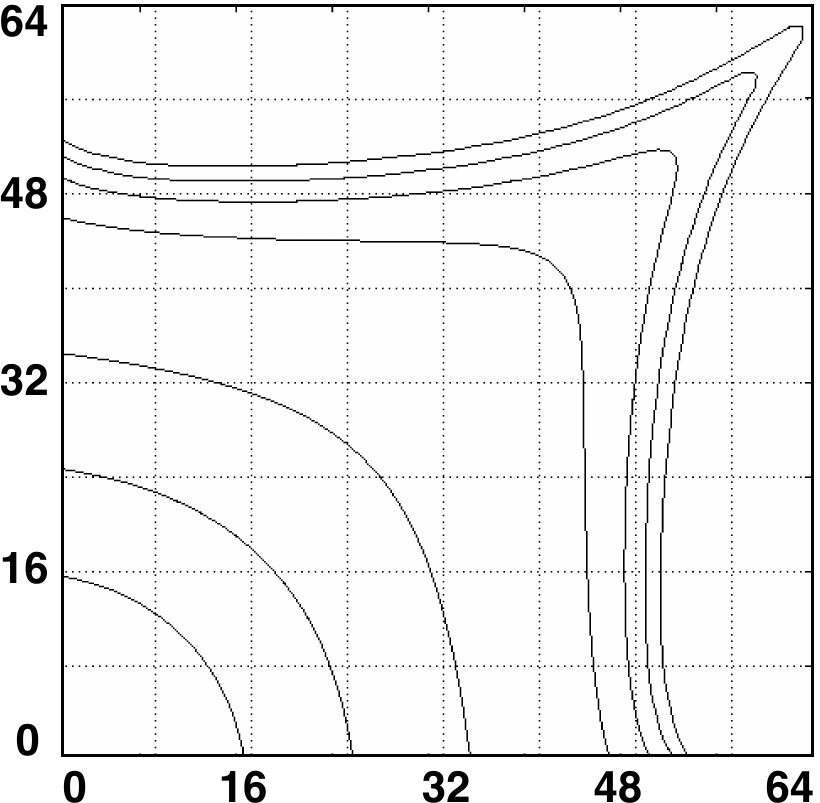}}
\\ 
\vspace{1.0cm}
\subfigure[\small NT2D: $128 \!\! \times \!\! 128$ células
]{\label{NT_dgr}
\includegraphics[scale=0.35]{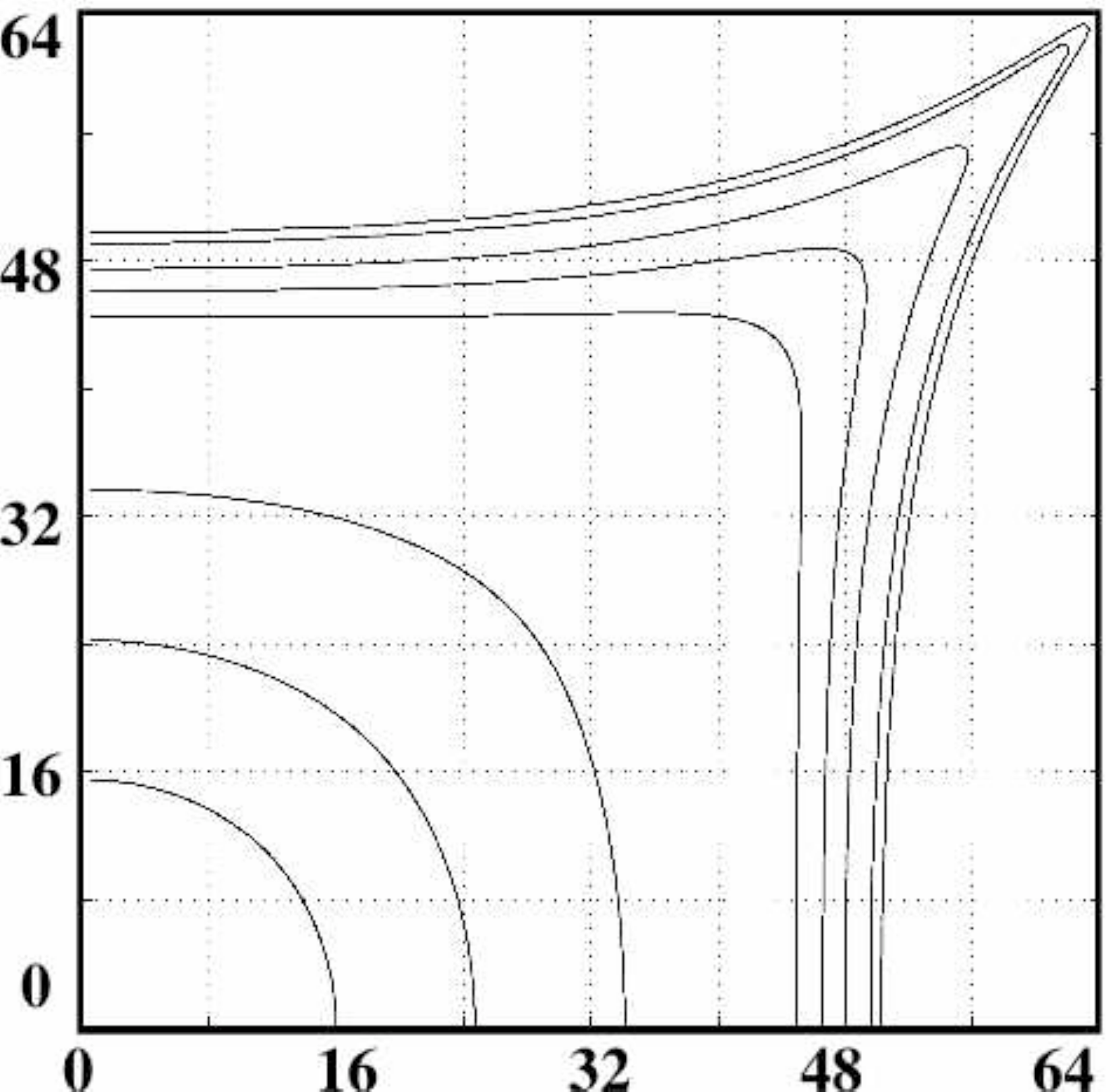}} 
\subfigure[\small KTdxd: $128 \! \times \! 128$ células
]{\label{NT_pgr}\quad\quad
\includegraphics[scale=0.7]{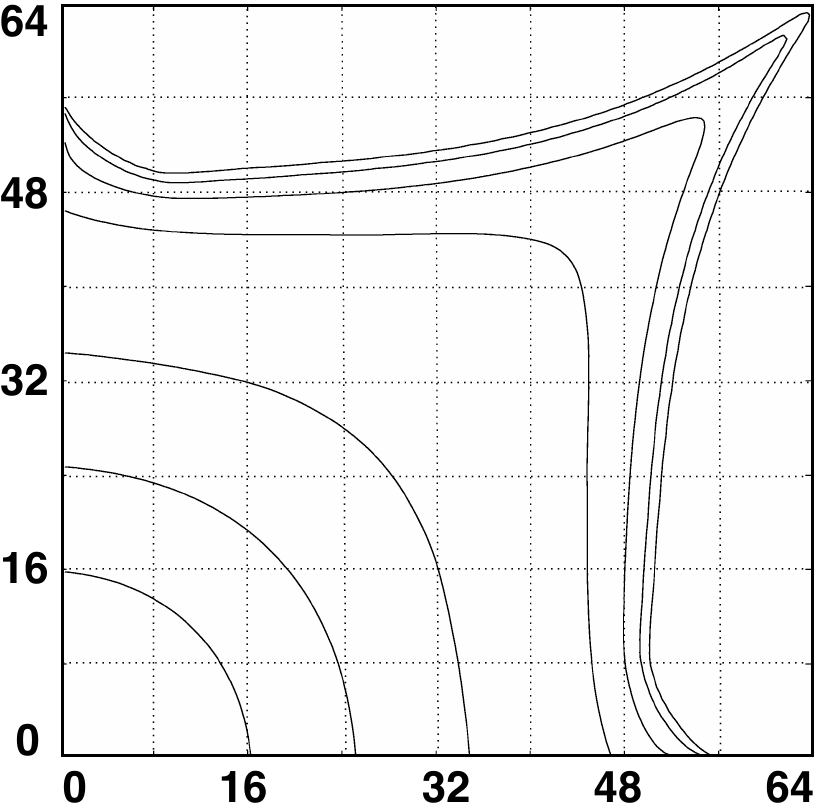}} \\
\caption{Water saturation level curves for two-phase flow in a five-spot well
configuration. The SD2-2D scheme was used in pictures 
in the left column and the KTdxd was used in pictures in the right column.}
\label{5_spot_nt_ktdxd}
\end{figure}

In the case of a five-spot flood with diagonal grid (Figure
\ref{5_spot_nt_ktdxd}), injection takes place at one corner and
production at the diametrically opposite corner; no flow is
allowed across the entirety of the boundary. In the case of
a five-spot flood with the parallel grid, 
injection takes place at two corners
(diametrically opposite), say left down and right up, and
production in the diametrically 'off corners', say right
down and left up.

\clearpage

\subsection{Analyzing the Numerical Results. Conclusions}

The Figures \ref{slab}, \ref{slab_CV1} and \ref{slab_CV2} refer to a
comparative study for the KT dimension by dimension and a
genuinely two-dimensional KT schemes showing the water
saturation surface plots after $350$ days of simulation for
three different values for the strength of the heterogeneity
of the fractal permeability field, $CV=0.5, 1.2$ and $2.2$.

Note that the genuinely two-dimensional KT scheme gives a
more accurate solution than the solutions computed by the KT
dimension by dimension scheme for the same grid. In fact we
observe that the KT dxd is only comparable in accuracy with
one degree of refinement (see Figures \ref{slab},
\ref{slab_CV1} and \ref{slab_CV2}).  The better accuracy of
the genuinely two-dimensional approach is due to a more
precise computation of the genuinely two-dimensional
numerical fluxes, with respect to the one dimensional
numerical fluxes in the dimension by dimension approach.

In the case of a 
five-spot geometry homogeneous reservoir, Figure
\ref{5_spot_nt_ktdxd} (diagonal grid)  shows the saturation level curves after $260$
days of simulation obtained with KTdxd and SD2-2D
schemes for two levels of spatial discretization.
In this figure \ref{5_spot_nt_ktdxd}, the
pictures on the left column are the results obtained with
the SD2-2D scheme and the ones on the right were computed
with the KTdxd scheme. In these Figures, the grid become
finer from top to bottom, having $64 \, \times \, 64$ and
$128 \, \times \, 128$ cells, respectively.

It is clear that the KTdxd scheme (right column pictures
in Figures \ref{5_spot_nt_ktdxd} is
producing incorrect boundary behavior. Moreover as the
computational grid is refined (right column and bottom
picture in Figure \ref{5_spot_nt_ktdxd}) 
this problem seems to be emphasized indicating that the
solution is not convergent.

The KTdxd scheme uses numerical fluxes in the $x$ and $y$
directions which can be viewed as generalizations of
one-dimensional numerical fluxes. We state that this type of
approximation for the fluxes leads to the incorrect boundary
behavior discussed above. This incorrect boundary behavior
led us to develop a new genuinely two-dimensional KT
scheme. The results obtained with this new scheme can be
seen in the left column of Figure \ref{5_spot_nt_ktdxd}. 
 It is clear that we have corrected the
boundary behavior just changing the approach from dimension
by dimension to our two-dimensional approach.  This fact indicates (but does not prove)
our idea that computing two-dimensional numerical fluxes
using straight generalizations of one-dimensional numerical
fluxes may produce incorrect numerical solutions.

\pagebreak


\def\cprime{$'$} \def\cprime{$'$}

\end{document}